\documentclass[a4paper,twoside,11pt]{amsart}
\newcommand{\Ueberschrift}{The Brauer--Manin obstruction for sections of the fundamental group} 
\newcommand{\Kurztitel}{Brauer--Manin for sections}
\usepackage[headings]{fullpage}
\usepackage{amsmath} \usepackage{amssymb} \usepackage{amsopn}
\usepackage{amsthm} \usepackage{amsfonts} \usepackage{mathrsfs}

\usepackage[pdfpagelabels,pdftex]{hyperref}
\usepackage[dvips]{graphicx} \usepackage[usenames]{color}
\usepackage[all]{xy}

\usepackage[OT2, T1]{fontenc}
\usepackage[russian,english]{babel}

\hypersetup{
  pdftitle={},
  pdfauthor={},
  pdfsubject={},
  pdfkeywords={},
  colorlinks=false,
  breaklinks=true,
  bookmarksopen=true,
  bookmarksnumbered=true,
  pdfpagemode=UseOutlines,
  plainpages=false}

\DeclareMathOperator{\rE}{E}
\DeclareMathOperator{\rF}{F}

\DeclareMathOperator{\rH}{H}
\DeclareMathOperator{\rI}{I}

\DeclareMathOperator{\rK}{K}

\DeclareMathOperator{\rN}{N}

\DeclareMathOperator{\rT}{T}

\DeclareMathOperator{\rc}{c}

\newcommand{\bA}{{\mathbb A}}

\newcommand{\bF}{{\mathbb F}}
\newcommand{\bG}{{\mathbb G}}

\newcommand{\bN}{{\mathbb N}}

\newcommand{\bP}{{\mathbb P}}
\newcommand{\bQ}{{\mathbb Q}}
\newcommand{\bR}{{\mathbb R}}

\newcommand{\bZ}{{\mathbb Z}}

\newcommand{\cD}{{\mathscr D}}

\newcommand{\cS}{{\mathscr S}}

\newcommand{\cU}{{\mathscr U}}

\newcommand{\cX}{{\mathscr X}}

\newcommand{\dO}{{\mathcal O}}

\newcommand{\fO}{{\mathfrak O}}

\newcommand{\fo}{{\mathfrak o}}
\newcommand{\fp}{{\mathfrak p}}

\DeclareSymbolFont{cyrletters}{OT2}{wncyr}{m}{n}
\DeclareMathSymbol{\Sha}{\mathalpha}{cyrletters}{"58}
\DeclareMathSymbol{\RussB}{\mathalpha}{cyrletters}{"42}

\newcommand{\surj}{\twoheadrightarrow} 
\newcommand{\inj}{\hookrightarrow}

\DeclareMathOperator{\id}{id}
\DeclareMathOperator{\pr}{pr}

\DeclareMathOperator{\Hom}{Hom}

\DeclareMathOperator{\Aut}{Aut}

\DeclareMathOperator{\coker}{coker}
\DeclareMathOperator{\im}{im}

\DeclareMathOperator{\gr}{gr}

\DeclareMathOperator{\Spec}{Spec}
\DeclareMathOperator{\Projj}{Proj}

\DeclareMathOperator{\divisor}{div}
\DeclareMathOperator{\Pic}{Pic}

\newcommand{\redu}{{\rm red}}

\newcommand{\OO}{\dO}

\DeclareMathOperator{\Alb}{Alb}
\DeclareMathOperator{\SAlb}{SAlb}

\newcommand{\Gm}{\bG_m}

\DeclareMathOperator{\ix}{index}
\DeclareMathOperator{\pe}{period}

\DeclareMathOperator{\Br}{Br}

\newcommand{\inv}{{\rm inv}}
\DeclareMathOperator{\Frob}{Frob}

\DeclareMathOperator{\res}{res}
\DeclareMathOperator{\cores}{cor}

\DeclareMathOperator{\pig}{\overline{\pi}}

\DeclareMathOperator{\Gal}{Gal}

\newcommand{\ep}{\varepsilon}
\newcommand{\ph}{\varphi}

\newcommand{\et}{\text{\rm \'et}}

\newcommand{\tame}{{\rm tame}} 

\newcommand{\sep}{{\rm sep}}
\newcommand{\alg}{{\rm alg}}
\newcommand{\nr}{{\rm nr}} 
\newcommand{\Sh}{{\rm sh}}

\newcommand{\tors}{{\rm tors}}
\newcommand{\ab}{{\rm ab}}

\newcommand{\bruch}[2]{\genfrac{}{}{0.5pt}{}{#1}{#2}}

\newcommand{\ov}[1]{\mbox{${\overline{#1}}$}} 

\newtheorem{thm}{Theorem}

\newtheorem{prop}[thm]{Proposition}
\newtheorem{lem}[thm]{Lemma}
\newtheorem*{lem*}{Lemma}

\newtheorem{cor}[thm]{Corollary}
\newtheorem{conj}[thm]{Conjecture}

\theoremstyle{definition}

\theoremstyle{remark}
\newtheorem{rmk}[thm]{Remark}

\newtheorem{ques}[thm]{Question}

\newenvironment{pro}[1][Proof]{{\it{#1:}} }{\hfill $\square$}
\newenvironment{pro*}[1][Proof]{{\it{#1:}} }{}
\newenvironment{pro**}[1][]{{\it{#1}} }{\hfill $\square$}

\newcounter{absatzcounter}[section]
\numberwithin{equation}{section}

\begin{document}

\hrule width\hsize
\hrule width\hsize
\hrule width\hsize
\hrule width\hsize

\vskip 0.2cm

arXiv version \hfill October 26, 2009

\vskip 0.2cm

\hrule width\hsize
\hrule width\hsize
\hrule width\hsize

\vspace{1.5cm}

\title[\Kurztitel]{\Ueberschrift} 
\author{ Jakob Stix}
\address{Jakob Stix, Mathematisches  Institut, Universit\"at Heidelberg, Im Neuenheimer Feld 288, 69120 Heidelberg}
\email{stix@mathi.uni-heidelberg.de}
\urladdr{http://www.mathi.uni-heidelberg.de/~stix/}
\thanks{The author acknowledges hospitality and support provided by the Isaac Newton Institute at Cambridge during the summer of 2009, when parts of this note were written.}

\subjclass[2000]{14H30, 11G20  Primary ; 14G05 Secondary }
\keywords{Section Conjecture, Rational points, Anabelian Geometry}
\date{October 26, 2009} 

\maketitle

\begin{quotation} 
  \noindent \small {\bf Abstract} --- We establish Grothendieck's section conjecture for an open subset of the Reichardt-Lind curve, and introduce the  notion of a Brauer--Manin obstruction for sections of the fundamental group extension.
\end{quotation}


\setcounter{tocdepth}{1} {\scriptsize \tableofcontents}


\section{Introduction} \label{sec:intro}

\noindent This note addresses the arithmetic of rational points on curves over  algebraic number fields with  the \'etale fundamental group as the principal tool.


\subsection{The section conjecture}
Let  $k$  be a field  with  absolute Galois group $\Gal_k = \Gal(k^\sep/k)$ with respect to a fixed separable closure $k^\sep$ of $k$. The fundamental group of a geometrically connected variety $X/k$ with geometric point $\ov{x} \in X$ sits in an extension
\begin{equation} \label{eq:pi1ext}
1 \to \pi_1(X \times_k k^\sep,\ov{x}) \to \pi_1(X,\ov{x}) \to \Gal_k \to 1,
\end{equation}
that  we abbreviate by $\pi_1(X/k)$ ignoring the basepoints right away in our notation.

To a rational point $a \in X(k)$ the functoriality of $\pi_1$ associates a splitting $s_a: \Gal_k \to \pi_1(X)$ of $\pi_1(X/k)$.  The functor $\pi_1$ depends a priori  on a pointed space which forces us to make choices and yields only a well defined $\pi_1(X \times_k k^\sep)$-conjugacy class  of sections.
The section conjecture of Grothendieck's contemplates about the converse.

\begin{conj}[see Grothendieck \cite{Groth:letter}] \label{conj:SC}
The map $a \mapsto s_a$ is a bijection of the set of rational points $X(k)$ with the set of 
$\pi_1(X \times_k k^\sep)$-conjugacy classes of sections of $\pi_1(X/k)$ if  $X$ is a smooth, projective, curve of genus  at least $2$ and $k/\bQ$ is finitely generated.
\end{conj} 
It was known to Grothendieck, that $a \mapsto s_a$ is injective by an application of the Mordell-Weil theorem in the number field case, see \cite{stix:cuspidalex} Appendix B, and \cite{Mz:localprop} Thm 19.1 for even a pro-$p$ result of injectivity. As we will deal in this note with problems of transfer from local to global, we emphasize, that Conjecture \ref{conj:SC} has also been conjectured for finite extensions $k/\bQ_p$. The argument for injectivity transfers almost literally to the $p$-adic case. 

It is also well known that in order to prove the section conjecture for curves over a fixed number field or $p$-adic field $k$ it suffices to show that having a section is the only obstruction for having a $k$-rational point, see \cite{koenigsmann} or \cite{stix:periodindex} Appendix C following a strategy of Nakamura and Tamagawa. In other words, if $\pi_1(X/k)$ splits, then we must guarantee a rational point. But we need this \textbf{weak version of the section conjecture} for all finite \'etale covers of $X$ that are geometrically connected over $k$.

A modified form of Conjecture \ref{conj:SC} for an open variety $U$ which is the complement of a divisor with normal crossings $Y$ in a geometrically connected proper smooth variety $X/k$  was proposed in \cite{Groth:letter}. The modification takes into account that $k$-rational points $y \in Y(k)$ contribute sections as follows. The natural surjection $D_y \surj \Gal_k$ of the decomposition group $D_y$ of $y$ in $\pi_1 U$ is known to admit splittings that lead by $D_y \subset \pi_1 U$ to sections of $\pi_1(U/k)$.


\subsection{Evidence for the section conjecture}
Recently there has been much work on the $p$-adic version of the conjecture and pieces of evidence for the section conjecture have emerged over the years. The most convincing evidence consists perhaps in Koenigsmann's proof in \cite{koenigsmann} of a birational analogue for function fields in one variable over a local $p$-adic field. 
A minimalist metabelian pro-$p$ approach to Koenigsmann's theorem was successfully developed by Pop in \cite{pop:propsc}, whereas Esnault and Wittenberg, \cite{ew} Prop 3.1, were able to reprove algebraically the "abelian part" via the cycle class of a section.

The analogue of the section conjecture over $\bR$, the real section conjecture, is known by  Mochizuki \cite{Mz:topics} Thm 3.13, following work of Sullivan, Cox and Huisman. Alternative proofs can be found  also in \cite{stix:periodindex} Appendix A, in the notably elegant simple argument of \cite{ew} Rmk 3.7(iv), using a fixed point theorem in \cite{pal} Thm 1.3, and moreover in  \cite{wickelgren:lcs} \S3.2, in which Wickelgren shows a geometrically 2-step nilpotent pro-$2$ version.

The first known examples of curves over number fields that satisfy the section conjecture are \textit{empty examples} in the sense that there are neither sections and hence nor rational points. The absence of sections is explained by the real section conjecture in case there are no real points on the curve for some real place. A local $p$-adic obstruction to sections is constructed in  \cite{stix:periodindex}, which explains the absence of sections in a number of cases when there are no $p$-adic points for some $p$-adic place.

Later Harari and Szamuely \cite{hs} gave examples of curves over number fields that are counterexamples to the Hasse principle, so they have local points everywhere but no global point, and nevertheless satisfy the section conjecture. These are again \textit{empty examples}, but as we mentioned above, the ostensibly dull case of empty curves is exactly the crucial class of examples.


\subsection{Outline of the paper} In this note, which contains parts of the author's Habilitations\-schrift \cite{stix:evidence} at the University Heidelberg,  we will discuss the analogue of the Brauer--Manin obstruction for an adelic point to be global, see Section \ref{sec:BMpoints}, for adelic sections of $\pi_1(X/k)$ and a number field $k$. This has been independently observed at least by O.~Wittenberg. Here we actually apply this obstruction to show in Section \ref{sec:RL} that the open subset $U$ of the Reichardt--Lindt curve described by 
\[
2y^2=x^4-17, \quad y \not=0
\]
does not admit a section of its fundamental group extension $\pi_1(U/\bQ)$. This case is particularly interesting further evidence for the section conjecture because $U$ has adelic points and the absence of sections is not explained by the strategy of Harari and Szamuely \cite{hs}. Furthermore, Section \ref{sec:RL} contains various generalizations: for arithmetic twists of the Reichardt--Lind curve and even for  an isotrivial family of such arithmetic twists inspired by an example by N.~D.~Elkies. The explicit examples for Brauer--Manin obstructions against sections are continued in Section \ref{sec:selmer} were  we illustrate the failure of the method to show the absence of sections in the case of the Selmer curve, which is given in homogeneous coordinates as 
\[
3X^3 + 4Y^3 + 5Z^3 = 0.
\]

In Section \ref{sec:prelim} we discuss the space of all sections $\cS_{\pi_1(X/k)}$. Among useful results for later parts of the note we generalize slightly as a byproduct the observation by P.~Deligne, see Appendix \ref{app:letterDeligne},  in that $\cS_{\pi_1(X/k)}$ is  naturally a compact pro-finite space for a proper variety $X$ over a number field $k$. 

Section \ref{sec:conditional} contains conditional results that depend on various conjectures. In particular we show that the local version of the section conjecture together with the conjecture that the Brauer--Manin obstruction is the only obstruction against rational points on curves imply the section conjecture for curves over number fields. 

\smallskip
\textbf{Acknowledgement}
The author is grateful to J.-L.~Colliot-Th\'el\`ene for explaining to him the Brauer--Manin obstruction to rational points on the Reichardt--Lind curve, and to N.~D.~Elkies for explaining his construction from the introduction of \cite{poonen:familyHasseP}. He thanks O.~Wittenberg for numerous comments and H.~Esnault for making available the letter from P.~Deligne to D.~Thakur, and is grateful to P.~Deligne for permitting the letter to appear as an appendix.


\section{Preliminaries on sections} \label{sec:prelim}
Let $X/k$ be a geometrically connected variety over a field $k$ of characteristic $0$ with algebraic closure $k^\alg$. The set of $\pi_1(X \times_k k^\alg)$-conjugacy classes of sections of $\pi_1(X/k)$ 
will be denoted by $\cS_{\pi_1(X/k)}$ and is called the \textbf{space of sections} of $X/k$.


\subsection{Base change}  For a field extension $K/k$ with algebraic closure $k^\alg \subseteq K^\alg$ the projection $X_K = X \times_k K \to X$ induces an isomorphism 
\[
\pi_1 (X_K \times_K K^\alg) \xrightarrow{\sim} \pi_1 (X \times_k k^\alg).
\]
Thus the extension $\pi_1(X_K/K)$ is the pullback of  $\pi_1(X/k)$ via the restriction $\Gal_K \to \Gal_k$ and 
\[
\pi_1(X_K) \xrightarrow{\sim} \pi_1 (X) \times_{\Gal_k} \Gal_K
\]
is an isomorphism. The defining property of a fibre product leads to a natural pullback map for spaces of sections 
\[
\cS_{\pi_1(X/k)} \to \cS_{\pi_1(X_K/K)}, \quad s \mapsto s_K = s \otimes K.
\]
We set $\cS_{\pi_1(X/k)}(K) = \cS_{\pi_1(X_K/K)}$ which describes a covariant functor on field extensions $K/k$ with values in sets.

\begin{lem} \label{lem:centraliser}
Let $k$ be either a number field or a  finite extension of $\bQ_p$, and let $X/k$ be smooth and geometrically connected. Any section $s : \Gal_k \to \pi_1(X)$ has an image with trivial centraliser in $\pi_1(X)$.
 \end{lem}
 \begin{pro} We argue by contradiction.
 As $\Gal_k$ is known to have trivial center, we may assume that $s(\Gal_k)$ centralises a nontrivial $\gamma \in \pi_1(X \times_k k^\alg)$. Let $H \lhd  \pi_1(X \times_k k^\alg)$ be a characteristic open subgroup with $\gamma \not\in H$. Then $\langle H,\gamma \rangle s(\Gal_k) \subseteq \pi_1(X)$ is an open subgroup corresponding to a finite \'etale cover $X' \to X$. The element $\gamma$ has nontrivial image in $\pi_1^\ab(X' \times_k k^\alg) = \langle H,\gamma \rangle^\ab$, as can be seen in the cyclic quotient $\langle H,\gamma \rangle/H$. Consequently, $\pi_1^\ab(X' \times_k k^\alg) $ contains a nontrivial $\Gal_k$-invariant subspace $\rH^0(k,\pi_1^\ab(X' \times_k k^\alg) )$. 
 
 Let $B$ be the generalized semiabelian Albanese variety of $X$. Then the Tate-module $\rT B = \varprojlim_n B[n]$ agrees as $\Gal_k$-module  with $\pi_1^\ab(X' \times_k k^\alg) $, so that 
 \[
 \rH^0(k,\pi_1^\ab(X' \times_k k^\alg) ) = \rH^0(k,\rT B) = \varprojlim_n B[n](k) = \varprojlim_n B(k)[n] = 0
 \]
 because the torsion group of $B(k)$ is finite if $k$ is a number field or a $p$-adic local field. This achieves a contradiction and completes the proof. 
\end{pro}

\begin{prop} \label{prop:galoisdescent}
Let $k$ be either a number field or a  finite extension of $\bQ_p$, and let $X/k$ be smooth and geometrically connected.

The map $E \mapsto \cS_{\pi_1(X/k)}(E)$ is a sheaf of sets on $\Spec(k)_\et$. In other words the space of sections satisfies Galois descent, i.e., for a Galois extension $E/F$ of finite extensions of $k$ the natural map 
\[
 \cS_{\pi_1(X/k)}(F) \to  \cS_{\pi_1(X/k)}(E)
\]
is injective and has the $\Gal(E/F)$-invariants as image.
\end{prop}
 \begin{pro}
 For the proof we may assume that $F=k$. We first show that restriction is injective. Let $s,t$ be classes of sections defined on $\Gal_k$ that agree when restricted on $\Gal_E$. We choose representatives of $s,t$ such that $s|_{\Gal_E}$ coincides with $t|_{\Gal_E}$. The difference
\[
a : \Gal_k \to \pi_1(X \times_k k^\alg), \qquad a_\sigma = t(\sigma)/s(\sigma)
\]
is a non-abelian cocyle, which means for all $\sigma,\tau \in \Gal_k$ we have
\[
a_{\sigma \tau} = a_\sigma s(\sigma)\big(a_\tau\big)s(\sigma)^{-1}.
\]
As $s$ and $t$ agree on $\Gal_E$, we have $a_\tau = 1$ for all $\tau \in \Gal_E$. It follows from putting $\tau \in \Gal_E$ and then $\sigma \in \Gal_E$ that $a$ factors over $\Gal(E/k)$ with values in the centraliser of $s(\Gal_E)$. By Lemma \ref{lem:centraliser} above we find $a \equiv 1$ and $s$ equals $t$ on all of $\Gal_k$.
 
The image of restriction is obviously contained in the set of $\Gal(E/k)$-invariant sections. Proving the converse, we assume that $s: \Gal_E \to \pi_1(X)$ is a Galois invariant sections. This means that for all $\sigma \in \Gal_k$ a suitable lift $\tilde{\sigma} \in \pi_1(X)$ will satisfy
\[
\big(\tilde{\sigma}(-)\tilde{\sigma}^{-1} \big)\circ s \circ  \big(\sigma^{-1}(-)\sigma\big) = s
\]
with the necessary conjugation by an element of $\pi_1(X \times_k k^\alg)$ being incorporated into the correct choice of $\tilde{\sigma}$. It follows that $\tilde{\sigma \tau}\big(\tilde{\sigma} \tilde{\tau}\big)^{-1}$
centralises $s(\Gal_E)$, hence equals $1$ by Lemma \ref{lem:centraliser}. For $\sigma \in \Gal_E$  we get
\[
s = \big(\tilde{\sigma}(-)\tilde{\sigma}^{-1} \big) \circ s \circ \big(\sigma^{-1}(-)\sigma\big) = \big(\tilde{\sigma}(-)\tilde{\sigma}^{-1} \big)  \circ \big(s(\sigma)^{-1}(-)s(\sigma) \big) \circ s
\]
so that $\tilde{\sigma} s(\sigma)^{-1}$ centralises $s(\Gal_E)$, hence $\tilde{\sigma} = s(\sigma)$ by yet another application of Lemma \ref{lem:centraliser}. It follows that $\sigma \mapsto \tilde{\sigma}$ is a section in $\cS_{\pi_1(X/k)}$ that extends the given section $s$ to all of $\Gal_k$ and this remained to be proved.
\end{pro}
 

\subsection{Reduction of sections} In this paragraph, let $S=\Spec(R)$ be the spectrum of an excellent  henselian discrete valuation ring with generic point $\eta = \Spec(k)$ and closed point $s=\Spec(\bF)$. For a geometrically connected, proper variety $U/k$, which is the generic fibre $j:U \subseteq X$ of a proper flat model $X/S$ with reduced special fibre $Y=X_{s,\redu}$, we have a specialisation map 
\[
U(k) = X(R) \to Y(\bF).
\]
The corresponding structure for sections is as follows.

\subsubsection{The ramification of a section}
The kernel of the natural map $\Gal_k \surj \pi_1(S,s) = \Gal_\bF$ is the \textbf{inertia group} $\rI_k$, which is the absolute Galois group of the  field of fractions  $k^\nr$ of the strict henselisation $R^\Sh$ with spectrum $\tilde{S} = \Spec(R^\Sh)$.  Inspection of the diagram
\[
\xymatrix@M+1ex@R-2ex{ 1 \ar[r] & \pi_1(U \times_k k^\alg) \ar[d] \ar[r] & \pi_1 (U) \ar@{->>}[d]^{\pi_1 (j)} \ar[r] & \Gal_k \ar@{->>}[d] \ar[r] & 1\\
1 \ar[r] & \pi_1(X \times_S \tilde{S}) \ar@{=}[d] \ar[r] & \pi_1 (X)  \ar@{=}[d]  \ar[r] & \pi_1(S,s)  \ar@{=}[d]  \ar[r] & 1\\
1 \ar[r] & \pi_1(Y \times_\bF \bF^\alg) \ar[r] & \pi_1 (Y) \ar[r] & \Gal_\bF \ar[r] & 1\\
}
\]
leads for each section $\sigma \in \cS_{\pi_1(X/k)}$ to  the \textbf{ramification homomorphism} 
\[
\rho = \rho_\sigma :\rI_k \to \pi_1(Y \times_\bF \bF^\alg),
\]
 that sits in a commutative diagram
\[
\xymatrix@M+1ex@R-2ex{ 1 \ar[r] & \rI_k \ar[d]^{\rho} \ar[r] & \Gal_k \ar[d]^{\pi_1 (j) \circ \sigma} \ar[r] & \Gal_{\bF} \ar[r] \ar@{=}[d]  & 1 \\
1 \ar[r] & \pi_1(Y \times_\bF \bF^\alg) \ar[r] & \pi_1 (Y) \ar[r] & \Gal_\bF \ar[r] & 1\\
}
\]
Consequently, the ramification $\rho$ is a $\Gal_\bF$-equivariant homomorphism. For sections $\sigma = s_u$ which belong to a rational point $u \in U(k)$ the ramification $\rho_{s_u}$ vanishes.

\subsubsection{No tame ramification} \label{sec:notameram}
The inertia $\rI_k$ is the semi-direct product of a pro-$p$ group, where $p$ is the residue characteristic, and the tame inertia  $\rI_k^\tame = \hat{\bZ}(1)\big(\bF^\alg\big) = \prod_{\ell \not= p} \bZ_\ell(1)$. By abuse of notion we will also call  tame inertia the images of splittings of $\rI_k \surj \rI_k^\tame$.

\begin{lem}  \label{lem:notameramification}
Let us assume in addition that $\bF$ is a finite field of characteristic $p$ and that $Y/\bF$ is a proper curve. Then $\im(\rho)$ is a pro-$p$ group, i.e.,  the restriction of $\rho$ to  tame inertia is trivial. 
\end{lem}
\begin{pro}
Otherwise, the image of tame inertia under $\rho$ would be infinite, because $\pi_1(Y \times_\bF \bF^\alg)$ has finite cohomological dimension. After replacing $Y$, and thus $U$ and $X$ by some finite \'etale cover, we would have a nontrivial $\Gal_\bF$-equivariant map 
\[
\rho|_{\rI_k^\tame} \ : \  \hat{\bZ}(1)\big(\bF^\alg\big) \to \pi_1^\ab(Y \times_\bF \bF^\alg) 
\]
that contradicts the Frobenius weights in \'etale cohomology of proper varieties over finite fields. Here $\pi_1^\ab$ is the abelianization of $\pi_1$ that is Pontrjagin dual to $\rH^1_\et$ and torsion free.
\end{pro}


\subsection{Topology on the space of sections} \label{sec:topology}

\subsubsection{Neighbourhoods} 
A \textbf{neighbourhood} of a section $s \in \cS_{\pi_1(X/k)}$ is an open subgroup $H$ of $\pi_1(X)$ together with a representative of $s$ whose image is contained in $H$ considered up to conjugation by $H \cap \pig_1(X)$. Equivalently, a neighbourhood is a finite \'etale map $h:X' \to X$ with $X'$ geometrically connected over $k$ and a section $s' \in \cS_{\pi_1(X'/k)}$ which maps to $s$ under $\pi_1(h)$. 

The images $U_{h} = \im\big(\pi_1(h): \cS_{\pi_1(X'/k)} \to  \cS_{\pi_1(X/k)}\big)$
where $h$ runs through the neighbourhoods of any section of $\pi_1(X/k)$ form the open sets of a topology on $\cS_{\pi_1(X/k)}$. We will consider $\cS_{\pi_1(X/k)}$ as a topological space endowed with this topology.

\subsubsection{Characteristic quotients} A topologically finitely generated pro-finite group $\Gamma$ is in several natural ways a projective limit of characteristic finite quotients along the index system $(\bN,<)$. For example, we may set $Q_n(\Gamma) = \Gamma/\bigcap_\ph \ker(\ph)$ where $\ph$ ranges over all continuous homomorphisms $\Gamma \to G$ with $G$ finite of order $\leq n$. Then $Q_n(\Gamma) $ is finite,  the natural map $\Gamma \to \varprojlim_n Q_n(\Gamma)$ is an isomorphism and the kernel of $\Gamma \surj Q_n(\Gamma)$ is preserved under any continuous automorphism of $\Gamma$.

As $\pi_1(X \times_k k^\alg)$ is topologically finitely generated if $k$ is of characteristic $0$ we might push $\pi_1(X/k)$ by the map $\pi_1(X \times_k k^\alg) \surj Q_n(\pi_1(X \times_k k^\alg))$ to obtain extensions $Q_n(\pi_1(X/k))$:
\[
1 \to Q_n(\pi_1(X \times_k k^\alg)) \to \pi_1(X)/\ker\left(\pi_1(X \times_k k^\alg) \surj Q_n(\pi_1(X \times_k k^\alg))\right)\to \Gal_k \to 1.
\]
We naturally extend the notation for the space of sections, so that $\cS_{Q_n(\pi_1(X/k))}$ denotes the $Q_n(\pi_1(X \times_k k^\alg))$-conjugacy classes of sections of $Q_n(\pi_1(X/k))$.
\begin{lem} \label{lem:quot}
The natural map
$\cS_{\pi_1(X/k)} \to \varprojlim_n \cS_{Q_n(\pi_1(X/k))}$
is a homeomorphism if we endow the right hand side with the pro-discrete topology.
\end{lem}
\begin{pro}
Let $s,t$ be representatives of clases in $\cS_{\pi_1(X/k)}$ that agree at every level $Q_n(\pi_1(X/k))$. Let $M_n \subset \pi_1(X \times_k k^\alg)$ be the nonempty set of elements which conjugate $s$ into $t$ at level $Q_n(\pi_1(X/k))$. Then $\varprojlim_n M_n$ is a projective limit of non-empty compact sets and is therefore non-empty. Any element in the limit conjugates $s$ into $t$. The surjectivity is clear.
\end{pro}

\subsubsection{The arithmetic case}
We now consider the case where $k$ is either a number field, or a $p$-adic local field or $\bR$.

\begin{prop}
Let $X$ be a  geometrically connected variety over either a $p$-adic or an achimedean local field $k$. Then the topological space of sections  $\cS_{\pi_1(X/k)}$ is a pro-finite set.
\end{prop}
\begin{pro} 
By Lemma \ref{lem:quot} it suffices to show finiteness of the set of sections of $Q_n(\pi_1(X/k))$, which follows in the local case from $\Gal_k$ being topologically finitely generated by \cite{jannsen} Thm 5.1(c), see also \cite{nsw} 7.4.1.
\end{pro}

\begin{prop} \label{prop:topology}
Let $X$ be a  proper, geometrically connected variety over an algebraic  number field $k$. Then the topological space of sections  $\cS_{\pi_1(X/k)}$ is a pro-finite set.
\end{prop}
\begin{pro}
By Lemma \ref{lem:quot} we need to show finiteness of the image of $\cS_{\pi_1(X/k)} \to  \cS_{Q_n(\pi_1(X/k))}$ for each $n$. Let $B_n \subset \Spec(\fo_k[\bruch{1}{n!}])$ be an open subset such that the extension $Q_n(\pi_1(X/k))$ is induced from an extension of $\pi_1(B_n)$ by $Q_n(\pi_1(X\times_k k^\alg))$.  Any section of $Q_n(\pi_1(X/k))$ coming from a section of $\pi_1(X/k)$ will be unramified at places in $B_n$ by Lemma \ref{lem:notameramification}. Any two such sections thus differ by a nonabelian cohomology class in $\rH^1(B_n,Q_n(\pi_1(X \times_k k^\alg))$, which is a finite set by Hermite's theorem that states the finiteness of the  set of field extensions of an algebraic number field with bounded degree and places of ramification.
\end{pro}


\begin{rmk}
Unfortunately, unlike with pro-finite groups, a countable pro-finite space need not be finite. In fact, the space $M=\{\bruch{1}{n} \ ; \ n \in \bN\} \cup \{0\}$ with the topology inherited as a subspace $M \subset \bR$ is a countable pro-finite set. Hence, contrary to a long term belief, the topological result of Proposition  \ref{prop:topology} does not reprove the Faltings--Mordell theorem of finiteness of the set of rational points on smooth, projective curves of genus at least $2$  once the section conjecture is known.
\end{rmk}

\begin{rmk} 
The results of Section \ref{sec:topology} have been obtained independently from the similar results by P.~ Deligne that are reported on in a letter from P.~Deligne to D.~Thakur. The author is grateful to P.~Deligne of having authorised its reproduction as an appendix for the convenience of the reader and for  giving due reference.
\end{rmk}


\subsection{Evaluation of units at sections} Let $U/k$ be a geometrically connected variety. 
For $n \in \bN$ not divisible by the characteristic, the Kummer sequence determines a homomorphism 
\[
\kappa: \OO^\ast(U) \to \rH^1(U,\mu_n) = \rH^1(\pi_1U,\mu_n), \quad f \mapsto \kappa_f.
\]
The \textbf{evaluation map} modulo $n$ for a section $s \in \cS_{\pi_1(U/k)}$ is the composite 
\[
\OO^\ast(U) \to k^\ast/(k^\ast)^n = \rH^1(k,\mu_n) = \rH^1(\Gal_k,\mu_n), \quad f \mapsto f(s) = s^\ast(\kappa_f).
\]
The functoriality of the Kummer sequence shows that for a section $s_u$ associated to a $k$-rational point $u \in U(k)$ we have
\[
f(s_u) = f(u) \mod (k^\ast)^n
\]
and so the name evaluation is justified. 

\subsubsection{Galois-equivariant evaluation}
Let $k'/k$ be a finite Galois extension with Galois group $G=\Gal(k'/k)$. For a section $s$ of $\pi_1(U/k)$ the base change $s' = s \otimes k'$ induces a $G$-equivariant map 
\[
s'^\ast \ : \  \rH^1(\pi_1 U_{k'},\mu_n) \to \rH^1(k',\mu_n)
\]
because for $\sigma \in  G$ we have 
\[
s'^\ast \circ \rH^1(\pi_1(\id \times \sigma)) = s'^\ast \circ \big(s(\sigma)^{-1}(-)s(\sigma) \big)^\ast = \big(s(\sigma)^{-1}s'(-)s(\sigma) \big)^\ast \]
\[ 
= \big(s' \circ (\sigma^{-1}(-)\sigma) \big)^\ast =  \big(\sigma^{-1}(-)\sigma\big)^\ast \circ s'^\ast  = \rH^1(\pi_1(\sigma)) \circ s'^\ast.
\]
\subsubsection{Evaluation and norms} \label{sec:evNorm}
More generally, for a finite field extension $k'/k$ the projection $\pr : U_{k'} \to U$ is finite flat and thus allows a norm map $N:\pr_\ast\Gm \to \Gm$. The induced map $\pr_\ast \mu_n \to \mu_n$ induces the corestriction on cohomology. Thus the following diagram commutes
\[
\xymatrix@M+1ex@R-2ex{\OO^\ast(U_{k'}) \ar[r]^(0.4)\kappa \ar[d]^N & \rH^1(\pi_1 U_{k'},\mu_n) \ar[r]^{s'^\ast} \ar[d]^\cores & \rH^1(k',\mu_n) \ar[d]^\cores \ar@{=}[r] & k'^\ast/(k'^\ast)^n \ar[d]^{N_{k'/k}} \\
\OO^\ast(U) \ar[r]^(0.4)\kappa  & \rH^1(\pi_1 U,\mu_n) \ar[r]^{s^\ast} & \rH^1(k,\mu_n) \ar@{=}[r] & k^\ast/(k^\ast)^n}
\]
so that $N(f)(s) = N_{k'/k}(f(s'))$ with $s' = s \otimes k'$ and $f \in \OO^\ast(U_{k'})$.


\section{Brauer--Manin obstructions} \label{sec:BM}


\subsection{Review of the Brauer-Manin obstruction for rational points} \label{sec:BMpoints}
The Brauer--Manin obstruction was introduced by Manin in \cite{manin:icm} and can explain the failure of the local--global principle. Let $k$ be an algebraic number field and $X/k$ a smooth,  geometrically connected variety. Let $k_v$ denote the completion of $k$ at a place $v$ so that the restricted product $\bA_k = \prod'_v k_v$ with respect to all places and the rings of intergers $\fo_v \subset k_v$ is the ring of ad\`els of $k$. 

On the set of adelic points $X(\bA_k)$ we introduce the equivalence relation with set of equivalence classes $X(\bA_k)_\bullet$, where two points are equivalent if they lie in the same connected component at each infinite place. 
A global cohomological Brauer class $A \in  \Br(X) \stackrel{\text{\tiny def}}{=} \rH^2(X,\Gm)$ yields a function
\[
\langle A, - \rangle : X(\bA_k)_\bullet \to \bQ/\bZ,
\quad 
x=(x_v)_v \mapsto \sum_v \inv_v(A(x_v))
\]
where $A(x_v)$ is the pullback of $A$ via $x_v$ to the Brauer group $\rH^2(k_v,\Gm)=\Br(k_v)$ that has the canonical invariant map $\inv_v:\Br(k_v) \to \bQ/\bZ$. The Hasse--Brauer--Noether local global principle for Brauer groups, see \cite{nsw} 8.1.17, i.e.,  the exactness of 
\begin{equation} \label{eq:eqHBN}
0 \to \Br(k) \to \bigoplus_v \Br(k_v) \xrightarrow{\sum_v \inv_v} \bQ/\bZ \to 0,
\end{equation}
shows that the global points $X(k)$ lie in the Brauer kernel 
\[
X(\bA_k)_\bullet^{\Br} := \{x \in X(\bA_k)_\bullet \ ; \ \langle A, x \rangle = 0 \text{ for all } A \in  \rH^2(X,\Gm) \}
\]
of all the functions $\langle A, - \rangle$. If a variety $X$ has adelic points but empty Brauer kernel $X(\bA_k)_\bullet^{\Br}$, then this explains the absence of global points $x \in X(k)$ and thus the failure of the local--global principle for $X$. For a more detailed exposition, see \cite{sko:torsors} \S5.2.

However, examples exist in dimension exceeding $1$, where the absence of rational points cannot be explained by a Brauer-Manin obstruction. The first was found by Skorobogatov \cite{sko} and could later be explained through a Brauer-Manin obstruction on a finite \'etale double cover. The recent examples of Poonen \cite{poonen} are shown to even resist being explained by the \'etale descent Brauer-Manin obstruction.


\subsection{Adelic sections}
We define the space of \textbf{adelic sections}  of $\pi_1(X/k)$ as the subset 
\[
\cS_{\pi_1(X/k)}(\bA_k) \subseteq \prod_v \cS_{\pi_1(X/k)}(k_v),
\]
with the product over all places of $k$, of all tuples $(s_v)$ such that for every homomorphism  $\ph: \pi_1 X \to G$ with $G$ finite the composites $\ph \circ s_v : \Gal_{k_v} \to G$ are unramified for almost all $v$. A group homomorphism $\Gal_{k_v} \to G$ is unramified if it kills the inertia subgroup $\rI_{k_v} \subset \Gal_{k_v}$.
\begin{prop}
The natural maps yield a commutative diagram
\[
\xymatrix@M+1ex@R-2ex{X(k) \ar[r] \ar[d] & X(\bA_k)_\bullet \ar[d]  \\
\cS_{\pi_1(X/k)} \ar[r] & \cS_{\pi_1(X/k)}(\bA_k) .}
\]
\end{prop}
\begin{pro}
For a section $s$ of $\pi_1(X/k)$ and a homomorphism $\ph: \pi_1 X \to G$ with $G$ finite, the composite $\ph \circ s: \Gal_k \to G$ describes a $G$-torsor over $k$ which therefore is unramified almost everywhere. Hence the tuple $(s \otimes k_v)_v$ is an adelic section.

Secondly, for an adelic point $(x_v) \in X(\bA_k)_\bullet$ the tuple of associated sections $s_{x_v}$ is adelic. Indeed, the $G$-torsor over $X$ corresponding to a homomorphism $\ph: \pi_1 X \to G$ with $G$ finite extends to a $G$-torsor over some flat model $\cX/B$ for a nonempty open $B \subset \Spec(\fO_k)$ for the integers $\fO_k$ in $k$. By the adelic condition for almost all places $v$ the point $x_v \in X(k_v)$ actually is a point $x_v : \Spec(\fo_v) \to  \cX$.
This means that $\ph \circ s_{x_v}$ factors as in the following diagram
\[
\xymatrix@M+1ex@R-2ex{
\Gal_{k_v} \ar[r]^(0.45){s_{x_v}}  \ar@{->>}[d] & \pi_1(X \times_k k_v) \ar@{}[r]|(0.55){\displaystyle \subseteq} \ar@{->>}[d] & \pi_1 (X) \ar[r]^{\ph} \ar@{->>}[d] & G \\
\pi_1(\Spec(\fo_v)) \ar[r]^{\pi_1(x_v)} & \pi_1(\cX \times_B \fo_v) \ar[r] & \pi_1 (\cX) \ar[ur] & 
}
\]
which shows that for such $v$ the composite $\ph \circ s_{x_v}$ is unramified.
\end{pro}
\begin{prop} \label{prop:alladelic}
Let $X/k$ be a proper geometrically connected variety. Then all tuples in $\prod_v \cS_{\pi_1(X/k)}(k_v)$ are adelic sections of $\pi_1(X/k)$.
\end{prop}
\begin{pro}
Let $\ph:\pi_1 X \to G$ be a homomorphism with $G$ finite and $(s_v)$ a tuple of local sections. Let $v$ be a place with residue characteristic $p \nmid \#G$, and such that the $G$-torsor corresponding to $\ph$ has good reduction over a flat proper model $\cX/B$ with $v \in B$. Then the restriction of $\ph \circ s_v$ to the inertia group $\rI_{k_v}$ factors over the ramification $\rho_{s_v}$ that by Lemma 
\ref{lem:notameramification} has a pro-$p$ group as its image and therefore dies in $G$. Hence for all such places $v$ the corresponding $\ph \circ s_v$ is unramified.
\end{pro}


\subsection{Brauer--Manin obstruction for sections}
Let $\alpha \in \rH^2(\pi_1X,\mu_n)$ be represented by the extension
\[
1 \to \mu_n \to E_\alpha \to \pi_1 X \to 1.
\]
The evaluation of $\alpha$ in a section $s_v \in \cS_{\pi_1(X/k)}(k_v)$ is the Brauer class $s_v^\ast(\alpha) \in \rH^2(k_v,\mu_n)$ that is represented by the pullback via $s_v$ of the extension $E_\alpha$.
\begin{prop} \label{prop:summable}
Let $(s_v)$ be an adelic section of $\pi_1(X/k)$. Then $s_v^\ast(\alpha)$ vanishes for all but finitely many places $v$.
\end{prop}
\begin{pro}
The extension $E_\alpha$ comes by inflation via a homomorphism $\ph : \pi_1 X \to G$ with $G$ finite from a class in $\rH^2(G,\mu_n)$ represented by an extension
\[
1 \to \mu_n \to E \to G \to 1.
\]
The class $s_v^\ast(\alpha)$ vanishes if and only if the composite $\ph \circ s_v$ lifts to $E$. The adelic condition on $(s_v)$ yields that for almost all $v$ the map $\ph \circ s_v$ factors over $\Gal_{\kappa(v)} \cong \hat{\bZ}$ which admits a lift to $E$ because $\hat{\bZ}$ is a free pro-finite group.
\end{pro}

Proposition \ref{prop:summable} allows the following definition of a function on adelic sections.
\[
\langle \alpha, - \rangle \ : \  \cS_{\pi_1(X/k)}(\bA_k)  \to \bQ/\bZ, \quad (s_v) \mapsto \sum_v \inv_v\big(s_v^\ast(\alpha)\big)
\]
\begin{thm} \label{thm:globalsectionsBMkernel}
The image of the natural map $\cS_{\pi_1(X/k)} \to \cS_{\pi_1(X/k)}(\bA_k)$ lies in the Brauer kernel
\[
 \cS_{\pi_1(X/k)}(\bA_k)^{\Br} := \{(s_v) \in   \cS_{\pi_1(X/k)}(\bA_k) \ ; \ \langle \alpha, (s_v) \rangle = 0 \text{ for all } n \in \bN \text{ and } \alpha \in  \rH^2(\pi_1X,\mu_n) \}.
\]
\end{thm}
\begin{pro}
For a section $s \in  \cS_{\pi_1(X/k)}$ and $\alpha \in \rH^2(\pi_1 X,\mu_n)$ we have 
\[
\langle \alpha, (s \otimes k_v)_v \rangle = \sum_v \inv_v\big((s \otimes k_v)^\ast(\alpha)\big) = \sum_v \inv_v\big(s^\ast(\alpha) \otimes_k k_v\big) = \big(\sum_v \inv_v\big)(s^\ast(\alpha)) = 0,
\]
that vanishes by the Hasse--Brauer--Noether local global principle for Brauer groups (\ref{eq:eqHBN}).
\end{pro}

The notation $ \cS_{\pi_1(X/k)}(\bA_k)^{\Br} $ requests a warning that its choice comes from the analogy with the Brauer kernel of adelic points and that in fact a more precise notation would refer to the obstruction coming from $\rH^2$ with coefficients in $\mu_n$. It it is by no means clear, albeit predicted by the local section conjecture, that the function $\langle \alpha, (s_v) \rangle$ on adelic sections only depends on the image $b(\alpha)$ in the Brauer group. But our choice of notation is shorter and suggestive and hopefully that sufficiently justifies its use.

The Kummer sequence yields an exact sequence
\begin{equation} \label{eq:Kummer}
0  \to \Pic(X)/n \Pic(X) \xrightarrow{\rc_1} \rH^2(X,\mu_n) \xrightarrow{b} \ _n\Br(X) \to 0
\end{equation}
where $\ _n\Br(X)$ is the $n$-torsion of $\Br(X)$. The composite $\rH^2(\pi_1 X,\mu_n) \subseteq \rH^2(X,\mu_n) \xrightarrow{b} \ _n\Br(X)$ is again denoted by $b$. The naturality of the Kummer sequence and $\rH^2(k_v,\mu_n) = \ _n\Br(k_v)$ imply 
\begin{equation} \label{eq:BMadelic}
\langle b(\alpha),(x_v)\rangle = \langle \alpha,(s_{x_v}) \rangle 
\end{equation}
and thus a commutative diagram.
\begin{equation} \label{eq:BMadelicdiagram}
\xymatrix@M+1ex@R-5ex{X(k) \ar[r] \ar[dd] & X(\bA_k)^{\Br}_\bullet \ar[dd]  \ar@{}[r]|{\displaystyle \subseteq} & X(\bA_k)_\bullet \ar[dd] \ar[dr]^(0.55){\langle b(\alpha),- \rangle} & \\
& & & \bQ/\bZ \\
\cS_{\pi_1(X/k)} \ar[r] & \cS_{\pi_1(X/k)}(\bA_k)^{\Br} \ar@{}[r]|{\displaystyle \subseteq} &  \cS_{\pi_1(X/k)}(\bA_k) \ar[ur]_(0.6){\langle \alpha,- \rangle}  & }
\end{equation}
\begin{prop} \label{prop:fibre2}
The middle facet  of diagram (\ref{eq:BMadelicdiagram}) is a fibre product for smooth, geometrically connected curves $X/k$.
\end{prop}
\begin{pro}
It suffices to prove that  $X(\bA_k)_\bullet $ is empty or that the map 
\[
b: \rH^2(\pi_1 X, \mu_n) \to \ _n\Br(X)
\]
 is surjective for every $n \geq 1$. Curves $X$ which are not a form of $\bP_k^1$ are algebraic $\rK(\pi,1)$ spaces and thus have $\rH^2(\pi_1 X,\mu_n) = \rH^2(X,\mu_n)$ so that the second condition follows from the Kummer sequence (\ref{eq:Kummer}).

A form of $\bP_k^1$ with local points everywhere is isomorphic to $\bP_k^1$ by the Hasse local--global principle. So the non-trivial forms $X$ have $X(\bA_k)_\bullet =\emptyset$. It remains to discuss $X=\bP_k^1$. But then the composite  
\[
\rH^2(k,\mu_n) = \rH^2(\pi_1 X,\mu_n) \subset \rH^2(X,\mu_n) \xrightarrow{b} \ _n\Br(X)
\]
is an isomorphism, because $\rc_1(\OO(1))$ kills the part of $\rH^2(X,\mu_n)$ which does not come from $\rH^2(\pi_1 X,\mu_n)$.
\end{pro}


\subsection{Locally constant Brauer--Manin obstructions}

The group of locally constant Brauer classes is defined as 
\[
\RussB(X) := \ker\Big(\Br(X) \to \bigoplus_v \Br(X \times_k k_v)/\Br(k_v)\Big).
\]
The evaluation of a Brauer class $A \in \RussB(X)$ on an adelic point does not depend on the choice of the local components and thus gives rise to a homomorphism 
\[
\iota_X : \ \RussB(X) \to \bQ/\bZ, \qquad A \mapsto \langle A, (x_v)\rangle
\]
independent of the choice of $(x_v) \in X(\bA_k)_\bullet$, see \cite{manin:icm}.  Let us define the group of $\RussB$-classes in $\rH^2(\pi_1 X,\bQ/\bZ(1))$ by 
\[
\rH^2_{\RussB}(\pi_1 X,\bQ/\bZ(1)) := \{\alpha \in \rH^2(\pi_1 X,\bQ/\bZ(1)) \ ; \ b(\alpha) \in \RussB(X) \}.
\]
The group $\rH^2_{\RussB}(\pi_1 X,\bQ/\bZ(1))$ contains the group of locally constant classes
\[
\rH^2_{\rm lc}(\pi_1 X,\bQ/\bZ(1)) := \ker\Big(\rH^2(\pi_1 X,\bQ/\bZ(1)) \to \prod_v  \rH^2(\pi_1 (X\times_k k_v),\bQ/\bZ(1)) /\rH^2(k_v,\bQ/\bZ(1)) \Big).
\]
The evaluation of a class $\alpha \in \rH^2_{\rm lc}(\pi_1 X,\bQ/\bZ(1))$ on an adelic section does not depend on the choice of the local components and thus gives rise to a homomorphism 
\[
j_X : \ \rH^2_{\rm lc}(\pi_1 X,\bQ/\bZ(1)) \to \bQ/\bZ, \qquad \alpha \mapsto \langle \alpha,(s_v)\rangle
\]
independent of the choice of the $(s_v) \in \cS_{\pi_1(X/k)}(\bA_k)$. If $X$ is a $\rK(\pi,1)$ with respect to cohomological degree $2$, more precisely if $b: \rH^2(\pi_1X,\bQ/\bZ(1)) \to \Br(X)$ is surjective, then the Brauer-Manin pairings yield a commutative diagram
\begin{equation} \label{eq:lcobst}
\xymatrix@M+1ex@R-2ex@C-2ex{X(\bA)_\bullet \ar[r] \ar@{^(->}[d] & \Hom(\RussB(X),\bQ/\bZ) \ar@{^(->}[d]  \ar[dr] & \\
\cS_{\pi_1(X/k)}(\bA_k) \ar[r] & \Hom(\rH^2_{\RussB}(\pi_1 X,\bQ/\bZ(1)), \bQ/\bZ) \ar[r]^{\res} & \Hom( \rH^2_{\rm lc}(\pi_1 X,\bQ/\bZ(1)),\bQ/\bZ)}
\end{equation}
with injective vertical maps. The image in the top row is $\iota_X$, while the image in the bottom row is $j_X$, regardless of the adelic point or section chosen. So if the image of the space of adelic sections in $\Hom(\rH^2_{\RussB}(\pi_1 X,\bQ/\bZ(1)), \bQ/\bZ)$ is not constant, then the local section conjecture fails, which gives a method of attack to falsify the local section conjecture.  On the other hand, if the restriction 
\[
b :  \ \rH^2_{\rm lc}(\pi_1 X,\bQ/\bZ(1)) \to \RussB(X)
\]
is still surjective, in which case we say that there are enough locally constant $\mu_n$-extensions, then the diagonal arrow is still injective. As $\iota_X$ gets mapped to $j_X$, and the latter vanishes if there is a global section $s \in \cS_{\pi_1(X/k)}$, under the assumption of having enough locally constant $\mu_n$-extensions we would find that the existence of a global section implies the vanishing of $\iota_X$. Consequently, an affirmative answer to the following question would be quite useful.
\begin{ques} \label{ques:locconst}
Does there exist enough locally constant $\mu_n$-extensions for curves?
\end{ques}


\subsection{Locally constant obstructions for curves} 

\begin{prop} \label{prop:sfdeg}
Let $X/k$ be a smooth, projective curve of genus $g \geq 1$ over a number field $k$, which admits adelic points. Then the Leray spectral sequence for $X \to \Spec(k)$ with coefficients in $\bQ/\bZ(1)$ degenerates at the $\rE_2$-level, at least for matters of $\rH^2$.
\end{prop}

\begin{rmk}
In fact, if $k$ has no real place, or restricted to the prime to $2$ component the whole spectral sequence of Proposition \ref{prop:sfdeg} degenerates at the $\rE_2$-level.
\end{rmk}

\begin{pro}
The cohomology $\rH^q(X \times_k k^\alg,\bQ/\bZ(1))$ vanishes for $q \geq 3$, and $\rH^3(k,\bQ/\bZ(1)) = 0$. Thus only the maps $d_2^{1,0}$ and $d_2^{2,0}$ have to be examined. The first is the Brauer obstruction map $\Pic_{X,\tors}(k) \to \Br(k)$ restricted to the torsion subgroup, see \cite{stix:periodindex} Appendix B. Hence its image is contained in the kernel of $\Br(k) \to \Br(X)$ which vanishes because of the local global principle for the Brauer group and the existence of local points which split the local maps $\Br(k_v) \to \Br(X \times_k k_v)$.
The map $d_2^{2,0}$ vanishes because the composite
\begin{equation} \label{eq:surF2}
\Pic(X) \otimes \bQ/\bZ \xrightarrow{\rc_1} \rH^2(X,\bQ/\bZ(1)) \to \rH^0(k,\rH^2(X \times_k k^\alg,\bQ/\bZ(1))) = \bQ/\bZ
\end{equation}
is surjective being nothing but the degree map tensored with $\bQ/\bZ$. 
\end{pro}

Here the map $\rc_1$ comes from the Kummer sequence on $X$ and sits in an exact seqeunce
\begin{equation} \label{eq:kummer}
0 \to \Pic(X) \otimes \bQ/\bZ \xrightarrow{\rc_1} \rH^2(X,\bQ/\bZ(1)) \to \Br(X) \to 0.
\end{equation}
Being a $\rK(\pi,1)$-space we may replace \'etale cohomology by group cohomology in the result of Proposition \ref{prop:sfdeg}. Let $\rF^\bullet\rH^2(\pi_1 X,\bQ/\bZ(1))$ be the filtration induced by the Leray (or in this case Hochschild--Serre) spectral sequence. We set
\[
\rF^1\rH^2_{\RussB}\big(\pi_1 X,\bQ/\bZ(1)\big) := \rF^1\rH^2\big(\pi_1 X,\bQ/\bZ(1)\big) \cap \rH^2_{\RussB}\big(\pi_1 X,\bQ/\bZ(1)\big).
\]


Let $k$ be an algebraic number field and let $A/k$ be a commutative group scheme. 
In order to fix a notation, we remind the definition of the Tate--Shafarevich group
\[
\Sha^1(k,A) := \ker\Big(\rH^1(k,A) \to \prod_v \rH^1(k_v,A)\Big).
\]
If $A$ is moreover divisible, e.g. for $A = \Pic_X^0$,
we define the Selmer group for  $A$ as the group
\[
\rH^1_{\rm Sel}(k,A) = \ker\Big(\rH^1(k,A_{\tors}) \to \prod_v \rH^1(k_v,A_{\tors})/ \delta_v\big(A(k_v) \otimes \bQ/\bZ\big)  \Big),
\]
where $A_{\tors}$ is the torsion subgroup, and 
 $\delta_v : A(k_v) \otimes \bQ/\bZ \to  \rH^1(k_v,A_{\tors})$ is the direct limit of connecting homomorphisms for the multiplication by $n$ sequence for $A$ on $X \times_k k_v$.
 

\begin{prop}  \label{prop:selmer}
Let $X/k$ be a smooth, projective curve of genus $g \geq 1$ over a number field $k$, which admits adelic points. Then the inclusion in $\rH^1(k,\Pic_{X,\tors})$ defines an isomorphism
\[
\rH^1_{\rm Sel}(k,\Pic^0_{X}) = \rF^1\rH^2_{\RussB}\big(\pi_1 X,\bQ/\bZ(1)\big) / \rH^2(k,\bQ/\bZ(1)).
\]
\end{prop}
\begin{pro}
The kernel $\tilde{P}$ of the map $\deg \otimes \bQ/\bZ : \Pic(X) \otimes \bQ/\bZ \to \bQ/\bZ $ sits in a short exact sequence
\[
0 \to \Pic_X^0(k) \otimes \bQ/\bZ \to \tilde{P} \to \bZ/\pe(X)\bZ \to 0
\]
From (\ref{eq:surF2}), (\ref{eq:kummer}) and using Tsen's theorem we deduce from 
\[
\xymatrix@M+1ex@R-2ex@C-3ex{
0 \ar[r] & \Pic_X^0(k) \otimes \bQ/\bZ \ar[d] \ar[r] & \rH^1_{\rm Sel}(k,\Pic^0_{X}) \ar@{^(->}[d]^i \ar[r] & \Sha^1(k,\Pic_X^0) \ar[r] \ar[d] & 0 \\
0 \ar[r] & \tilde{P} 
\ar[r] & \rF^1\rH^2_{\RussB}(\pi_1 X,\bQ/\bZ(1))/\rH^2(k,\bQ/\bZ(1)) \ar[r] & \Sha^1(k,\Pic_X) \ar[r] & 0
}
\]
via the snake lemma an exact sequence
\[
0 \to \bZ/\pe(X)\bZ \xrightarrow{\delta} \bZ/\pe(X)\bZ \to \coker(i) \to 0.
\]
This proves the claim.
\end{pro}

As a consequence of Proposition \ref{prop:selmer}, we obtain the following variant of diagram (\ref{eq:lcobst}) 
\begin{equation}
\xymatrix@M+1ex@R-2ex{X(\bA)_\bullet \ar[r] \ar@{^(->}[d] & \Hom(\RussB(X)/\Br(k),\bQ/\bZ) \ar@{^(->}[d]  \\
\cS_{\pi_1(X/k)}(\bA_k) \ar[r] & \Hom(\rH^1_{\rm Sel}(k,\Pic^0_{X}) , \bQ/\bZ) }
\end{equation}
where the vertical maps are still injective, and the top horizontal arrow is constant with value the induced map by $\iota_X$. 

Following Sa\"\i di and Tamagawa we call a local section $s_v$ \textbf{good} if the composite
\[
\Pic^0_X(k_v) \otimes \bQ/\bZ \to \rH^2(\pi_1 X \times_k k_v, \bQ/\bZ(1)) \xrightarrow{s_v^\ast} \rH^2(k_v,\bQ/\bZ(1))
\]
vanishes. The subset $\cS^{\rm good}_{\pi_1(X/k)}(\bA_k)$ of the adelic sections such that all local components are good in the sense above receives the adelic sections corresponding to adelic points. Moreover, all good adelic section induce the same homomorphism $\rH^1_{\rm Sel}(k,\Pic^0_{X}) \to  \bQ/\bZ$ as their locally constant Brauer obstruction, which therefore equals 
\[
\iota_X : \ \rH^1_{\rm Sel}(k,\Pic^0_{X}) \surj \RussB(X)/\Br(k) \to  \bQ/\bZ.
\]
We conclude that if the adelic section associated to a section $s$ of $\pi_1(X/k)$ is good, then $\iota_X$ vanishes by Theorem \ref{thm:globalsectionsBMkernel}. The prime to $p$-part of a local section is good if the residue characteristic is $p$. The goodness of the $p$-part however remains largely a mystery.

\begin{rmk}
The natural map describes an inclusion 
\[
 \rH^2_{\rm lc}(\pi_1 X,\bQ/\bZ(1))/\rH^2(k,\bQ/\bZ(1)) \subset  \rH^1_{\rm Sel}(k,\Pic^0_{X})
\]
which, when analysed carefully, makes a positive answer to Question \ref{ques:locconst} very unlikely.
\end{rmk}


\section{Conditional results}  \label{sec:conditional}


\subsection{The assumptions}

In Section \ref{sec:conditional}, and in this section only, we work under the following hypotheses when neccessary.
\begin{itemize}
\item[(BM)] The Brauer--Manin obstruction is the only obstruction that prevents the existence of rational points  on smooth, projective geometrically connected curves  over algebraic number fields.
\item[(T$\Sha$)] The Tate-Shafarevic group $\Sha^1(k,\Pic_X^0)$ of the jacobian $\Pic_X^0$ of a smooth projective curve $X$ over a number field $k$ is finite.
\item[(LSC)] The local section conjecture holds true, i.e., the section conjecture for a base field that is finite over some $\bQ_p$. 
\end{itemize}

The credibility of the above assumptions varies. Only recently, several authors contributed evidence for (BM), see \cite{scharaschkin:brm},  \cite{flynn}, \cite{poonen:heuristics}, \cite{poonenvoloch} and \cite{stoll:finitedescent}. For surfaces and beyond (BM) was expected to fail a long time ago, and an example was finally provided by Skorobogatov, see \cite{sko} and \cite{poonen}. The assumption (T$\Sha$) belongs to main stream mathematics for quite some time and enters for example in the Birch and Swinnerton-Dyer conjecture or even more intrinsically in deeper conjecture about motives. The last assumption (LSC) is of course a necessity in light of any progress for the section conjecture over number fields which progresses along the avenue of the interaction between local and global in number theory. As weak evidence towards (LSC) we might regard the recent works  \cite{koenigsmann}, \cite{pop:propsc} and \cite{stix:periodindex}.


\subsection{\texorpdfstring{$0$}{0}-cycles instead of rational points} 
The Brauer--Manin obstruction against rational points can be extended to an obstruction against $k$-rational $0$-cycles. 

\subsubsection{Cassels--Tate pairing and locally constant obstructions} \label{sec:CTlco}
Let $X/k$ be a smooth, projective geometrically connected curve. The  Leray spectral sequence for $X \to \Spec(k)$ yields an exact sequence 
\[
\Br(k) \to \RussB(X) \to \Sha^1(k,\Pic_X) \to 0.
\]

\begin{lem} \label{lem:lcbrauerandCT}
Let $X/k$ be a smooth projective geometrically connected curve with adelic points and let $X \to W = \Alb_X^1$ be the map into its universal torsor under an abelian variety. There is a commutative diagram
\[
\xymatrix@R+1ex@R-2ex{\RussB(X) \ar@{->>}[r] \ar[dr]_{\iota_X} & \Sha^1(k,\Pic_X) \ar[d]^\delta & \ar@{->>}[l] \Sha^1(k,\Pic_X^0) \ar[dl]^{\langle [W],- \rangle_{\rm CT}}\\
 & \bQ/\bZ & }
\]
where $\langle [W],- \rangle_{\rm CT}$ is the Cassels--Tate pairing 
\[
 \Sha^1(k,\Alb_X) \times \Sha^1(k,\Pic^0_X) \to \bQ/\bZ
\]
evaluated at the class $[W] \in  \Sha^1(k,\Alb_X)$.
\end{lem}

\begin{pro}
This can be found for example in \cite{manin:icm} Theorem 6, Proposition 8c, \cite{milne:BrSha} Lemma 2.11, or \cite{sko:torsors} Theorem 6.2.3.
\end{pro}

Because the Cassels--Tate pairing is non-degenerate modulo the subgroup of divisible elements we get the following easy corollary, see D.~Eriksson and V.~Scharaschkin \cite{erikssonscharaschkin} Thm 1.2, who were the first to state that locally constant classes in the Brauer group suffice to detect the absence of $0$-cycles of degree $1$ on curves over number fields with finite Tate-Shafarevich groups.
 
\begin{cor}[Eriksson--Scharaschkin \cite{erikssonscharaschkin}] \label{cor:index1}
Under the assumptions of Lemma \ref{lem:lcbrauerandCT}, if (T$\Sha$) holds for $\Pic_X^0$, then the following are equivalent.
\begin{itemize}
\item[(a)] $\iota_X = 0$,
\item[(b)] $\Alb^1_X$ has a rational point,
\item[(c)] ${\rm period}(X) = 1$,
\item[(d)]  ${\rm index}(X) = 1$.
\end{itemize}
\end{cor}
\begin{pro}
If $X$ has an adelic point, then the kernel of $\Br(k) \to \Br(X)$ is trivial but surjects onto $\pe(X)\bZ/\ix(X)\bZ$ by the degree map. Hence (c) is equivalent to (d). By definition (c) is equivalent to (b) which under the assumtption (T$\Sha$) is equivalent to the triviality of the map $\langle [W], - \rangle_{\rm CT}$, hence equivalent to $\iota_X=0$ by Lemma \ref{lem:lcbrauerandCT}.
\end{pro}


\subsubsection{Locally constant obstructions assuming finiteness of \texorpdfstring{$\Sha$}{Sha}}

\begin{thm}
Let $X/k$ be a smooth, projective geometrically connected curve of genus at least $1$ with adelic points and such that $\pi_1(X/k)$ splits. Under the assumption (LCS) that the local section conjecture is true  for $X/k$ and (T$\Sha$) that the Tate-Shafarevich group for $\Pic_X^0$ is finite, we find that $X$ has index $1$, i.e., there is a rational $0$-cycle of degree $1$ on $X$.
\end{thm}
\begin{pro}
Under the assumption of the local section conjecture, in diagram (\ref{eq:lcobst}) the left vertical map is a bijection, showing that the image of $\iota_X$ in $ \Hom(\rH^2_{\RussB}(\pi_1 X,\bQ/\bZ(1)), \bQ/\bZ)$ equals the image of the adelic section $(s \otimes k_v)$ coming from the global section $s \in \cS_{\pi(X/k)}$. Hence $\iota_X$ vanishes by Theorem \ref{thm:globalsectionsBMkernel}. By Corollary \ref{cor:index1} we may deduce that $X$ has index $1$.
\end{pro}


\subsection{A fibre square and a local--global principle}

\begin{prop} \label{prop:BMkernelsequal}
Under the assumption (LSC) that the local section conjecture holds for smooth, projective curves over finite extensions of $\bQ_p$ of genus at least $2$ and all $p$, then the natural map 
\[
X(\bA_k)^{\Br}_\bullet  \to \cS_{\pi_1(X/k)}(\bA_k)^{\Br}
\]
is a bijection for a smooth, projective geometrically connected curve $X/k$ of genus at least $2$ over an algebraic number field $k$.
\end{prop}
\begin{pro} 
The assumption of the local section conjecture and the validity of the real section conjecture, see \cite{stix:periodindex} Appendix A, yield a bijection $X(\bA_k)_\bullet  \to \cS_{\pi_1(X/k)}(\bA_k)$. 
The result then immediately follows from Proposition \ref{prop:fibre2}.
\end{pro}

\begin{prop} \label{prop:fibre1}
Let $X$ be smooth, projective, geometrically connected curve over the number field $k$.
Under the assumption (BM) that the Brauer--Manin obstruction is the only obstruction that prevents the existence of rational points  on smooth, projective geometrically connected curves  over algebraic number fields, the left facet of diagram (\ref{eq:BMadelicdiagram})
\begin{equation} \label{eq:fibreproduct}
\xymatrix@M+1ex@R-2ex{X(k) \ar[r] \ar[d] & X(\bA_k)^{\Br}_\bullet \ar[d]   \\
\cS_{\pi_1(X/k)} \ar[r] & \cS_{\pi_1(X/k)}(\bA_k)^{\Br} }
\end{equation}
is a fibre product square.
\end{prop}
\begin{pro}
The vertical maps are injective by the known injectivity part of the global and the local section conjecture.

It suffices to show that a section $s$ of $\pi_1(X/k)$ that locally comes from an adelic point $(x_v)$ forces the existence of a rational point $x \in X(k)$. Indeed,  we may look at neighbourhoods $h: X' \to X$ of $s$ with a lift $t\in \cS_{\pi_1(X'/k)}$ of $s$. Because of  $s \otimes k_v = s_{x_v}$ the lift $t$ determines a unique lift of the adelic point $(x_v)$ to an adelic point of $X'$ that moreover is in $X'(\bA_k)_\bullet^{\Br}$ by Theorem \ref{thm:globalsectionsBMkernel} and Proposition \ref{prop:fibre2}.

It follows that the existence of rational points on $X'$ is not Brauer--Manin obstructed whereby we find a rational point $x' \in X'(k)$ by assumption. The usual limit argument in the tower of all neighbourhoods shows that the sections $s_{h(x')}$ converge to $s$ and that  the limit of the $h(x')$ becomes a rational point $x \in X(k)$ with $s_x=s$, see  \cite{stix:periodindex} Appendix C.
\end{pro}


The following logical dependence between different conjectures has been independently observed by O.~Wittenberg.
\begin{thm} \label{thm:LSCBMandSC}
If the local section conjecture (LSC) holds true for smooth, projective curves over finite extensions of $\bQ_p$ of genus at least $2$ and all $p$, and (BM) if the Brauer--Manin obstruction is the only obstruction that prevents the existence of rational points over algebraic number fields on smooth, projective curves, then the prediction of the section conjecture holds true: for a smooth, projective geometrically connected curve $X/k$ of genus at least $2$ over an algebraic number field $k$ the natural map
\[
X(k) \to \cS_{\pi_1(X/k)}
\]
is a bijection.
\end{thm}
\begin{pro} 
This is an immediate consequence of Proposition \ref{prop:BMkernelsequal} and Proposition \ref{prop:fibre1}.
\end{pro}


For a discussion of the finite descent obstruction,  the close cousin of the Brauer--Manin obstruction, on adelic points instead of sections and an analogue of Theorem \ref{thm:LSCBMandSC} we refer to the work of M.~Stoll \cite{stoll:finitedescentdraft8} \S9.


\section{Application: the Reichardt--Lind curve} \label{sec:RL}

\subsection{Geometry and arithmetic of the Reichardt--Lind curve}  \label{sec:geoarithRL}
For the convenience of the reader we recall the geometry and arithmetic of the Reichardt--Lind curve. 
For an elementary exposition on the Reichardt--Lind curve we may refer to \cite{aitkenlemmermeyer}. 
The \textbf{affine Reichardt--Lind curve} $U$ is the affine curve in $\bA_\bQ^2$ given by the equation
\begin{equation} \label{eq:affRL}
2Y^2 = Z^4-17,\quad Y \not= 0.
\end{equation}
The \textbf{projective Reichardt--Lind curve} $X$ is the projective curve over $\bQ$ given by the homogeneous equations in $\bP_\bQ^3$
\[
\left\{\begin{array}{l}
2T^2 = A^2 - 17 B^2 \\
AB = C^2
\end{array}\right.
\]
The map 
\[
(y,z) \mapsto [a:b:c:t]=[z^2:1:z:y]
\]
identifies $U$ with the open subset of $X$ defined by $TB \not= 0$. The curve $X$ has a model 
\[
\cX=\Projj \left(\bZ[\bruch{1}{34}][A,B,C,T]/(2T^2=A^2-17B^2,\ AB=C^2)\right)
\]
 that is smooth, projective with geometrically connected fibres above $\Spec \bZ[\bruch{1}{34}]$ as can be checked via the jacobian criterion: the matrix
\[
\left(\begin{array}{rrrr}
2A & -34B & 0 & -4T \\
B & A & -2C & 0
\end{array}\right)
\]
has full rank everywhere on $\cX$. We have $\Omega_X = \OO(-4+2+2)|_X = \OO_X$ by the adjunction, so that $X$ is a smooth, projective curve of genus $1$ over $\bQ$ with good reduction outside $\{2,17\}$. 
The divisor $\cD \subset \cX$ given by the reduced locus of $TB=0$ is isomorphic to 
\[
\cD=\Spec (\bZ[\bruch{1}{34},\sqrt[4]{17}]) \amalg \Spec( \bZ[\bruch{1}{34},\sqrt{2}] )
\]
and is relative effective and \'etale over $\Spec(\bZ[\bruch{1}{34}])$ of degree $6$. The complement $\cU = \cX -\cD$ is a flat model of $U$ over $\Spec(\bZ[\bruch{1}{34}])$ as a hyperbolic curve of type $(1,6)$. We have 
\[
\cU = \Spec \left(\bZ[\bruch{1}{34}][Y^{\pm 1},Z]/(2Y^2=Z^4-17)\right)
\]
and $\cU$ is precisely the locus where $Y \in \OO_\cU$ is invertible as a rational function on $\cX$. The divisor $\divisor(Y)$ equals $P-2\cdot Q$ with $P \cong \Spec (\bZ[\bruch{1}{34},\sqrt[4]{17}])$ and $Q \cong  \Spec( \bZ[\bruch{1}{34},\sqrt{2}])$.

The projective Reichardt--Lind curve $X$ is a principal homogeneous space under its Jacobian $E=\Pic_X^0$. The equation for $E/\bQ$ can be found in \cite{macrae} p.~278, namely $E$ is given by the affine equation
\[
Y^2 = X(X^2+17)
\]
The elliptic curve $E/\bQ$ is encoded by $18496k1$ with coefficient vector $[0,0,0,17,0]$ on Cremona's list, see \cite{cremona}, and has conductor $N=18496=2^6\cdot 17^2$, $E(\bQ) = \bZ/2\bZ$, analytic rank $0$, hence finite $\Sha^1(\bQ,E)$ of analytic order $4$.

The significance of the  Reichardt--Lind curve $X$ lies in the fact that  its class in $\rH^1(\bQ,E)$ belongs to the group $\Sha^1(\bQ,E)$ of $E$-torsors that are locally trivial at every place of $\bQ$. This is clear outside $2,17$ and $\infty$ by good reduction, at $\infty$ it is obvious, at $17$ we have $Q \in X(\bQ_{17})$, and at $2$  we have $P \in X(\bQ_2)$. That $X \in \Sha^1(\bQ,E)$ represents a non-trivial class means that $X(\bQ) = \emptyset$  and is recalled in Section \ref{sec:absenceRP}.
As such, the curve $X$ violates the Hasse local--global principle and was discovered  in 1940 by Lind \cite{lind}, and independently  in 1942 by  Reichardt \cite{reichardt}.


\subsection{Review of the classical proof for absence of rational points} \label{sec:absenceRP}

\subsubsection{Via Gau\ss\ reciprocity}
A solution in $\bQ$ of $2Y^2 = Z^4-17$ leads after clearing denominators to the equation  $2y^2 = z_0^4 - 17z_1^4$ with integers $z_0,z_1,y$ and $z_0,z_1$ coprime. Then $17 \nmid y$ and thus for each odd prime $p \mid y$ we find that $17$ is a square modulo $p$, hence $\Big(\bruch{p}{17}\Big) = \Big(\bruch{17}{p}\Big) = 1$. As also $\Big(\bruch{-1}{17}\Big)=1$ and $\Big(\bruch{2}{17}\Big)=1$, we find that $y$ is a square modulo $17$. But then $2$ must be a $4^{th}$ power modulo $17$, which it is not by $2^4 \equiv -1 \mod 17$, contradiction.

\subsubsection{Via a Brauer--Manin obstruction} \label{sec:absenceRPBM} 
The absence of global rational points  on $X$ can also be explained by a Brauer--Manin obstruction. The Kummer sequence with multiplication by $2$ yields a map 
\[
\OO^\ast(\cU)/(\OO^\ast(\cU))^2 \to \rH^1(\cU,\mu_2), \quad t \mapsto \chi_t
\]
and we use the cup-product $\alpha_{\cU} = \chi_Y \cup \chi_{17} \in \rH^2(\cU,\mu_2^{\otimes 2})$. The localization sequence together with the Gysin isomorphism yields an exact sequence
\[
 \rH^2(\cX,\mu_2^{\otimes 2}) \to \rH^2(\cU,\mu_2^{\otimes 2}) \xrightarrow{\res} \rH^1(\cD, \mu_2).
\]
The residue $\res(\alpha_\cU)$ can be computed with the help of excision or via the tame symbol of Milnor K-theory. Near $P \cong \Spec (\bZ[\bruch{1}{34},\sqrt[4]{17}])$ the character $\chi_{17}$ vanishes as $17$ becomes a square, and near $Q \cong  \Spec( \bZ[\bruch{1}{34},\sqrt{2}])$ the function $Y$ becomes a square and so $\chi_Y$ vanishes. It follows that $\alpha_\cU$ lifts to an element $\alpha \in \rH^2(\cX,\mu_2)$, where we have identified $\mu_2^{\otimes 2}$ with $\mu_2$.

We will compute the function $\langle A, - \rangle$ on adelic points $X(\bA_{\bQ})_\bullet$ for the image $A$ of $\alpha$ under $\rH^2(\cX,\mu_2) \to \Br(\cX)$. Let $(x_v) \in X(\bA_{\bQ})_\bullet$ be an adelic point. Outside $2,17$ and $\infty$ the class $A$ has good reduction and thus $A(x_v)$ belongs to $\Br(\fo_v)=0$. At $2$ and $\infty$ we find that $17$ is a square and because by continuity of $\inv_v(A(x_v))$ we may assume that actually $x_v \in U(\bQ_v)$ we find for $v=2,\infty$ that 
\[
\inv_v(A(x_v)) = \inv_v(x_v^\ast(\alpha)) = \inv_v(x_v^\ast(\chi_Y \cup \chi_{17})) = 0.
\]
It remains to compute $\langle A, (x_v) \rangle = \inv_{17}(A(x_{17}))$ where $A(x_{17})$ is the quaternion algebra $(y,17) = \Big(\bruch{y,17}{\bQ_{17}}\Big)$ with $y=Y(x_{17})$.  We set $E=\bQ_{17}(\sqrt[4]{17})$ and $F=\bQ_{17}(\sqrt{17})$. The norm residue property states that  under the identification 
\[
\bruch{1}{2}\bZ/\bZ = \Gal(F/\bQ_{17}) = \bQ_{17}^\ast/N_{F/\bQ_{17} } \big(F^\ast\big).
\]
we find $\inv_{17}((y,17)) = y \mod N_{F/\bQ_{17} } \big(F^\ast\big)$.
The field extension $E/\bQ_{17}$ is abelian and totally ramified of degree $4$. 
Multipliction by $2$ therefore induces a map 
\[
\bQ_{17}^\ast/N_{F/\bQ_{17}}  \big(F^\ast\big) \to \bQ_{17}^\ast/N_{E/\bQ_{17}}  \big(E^\ast\big)
\]
that by local class field theory is a map $\mu_2 \to \mu_4$ with image the squares in $\mu_4$, hence it is injective. As $y$ satisfies $2y^2 = z^4 - 17$ for $z=Z(x_{17})$ we see that $2y^2$ is a norm from $E$. We conclude that the image of $\inv_{17}((y,17))$ in $ \bQ_{17}^\ast/N_{E/\bQ_{17}}  \big(E^\ast\big)$ is given by the class of $1/2$ which is nontrivial because $2$ is not a $4^{th}$ power in $\bF_{17}^\ast$ and by local class field theory
\[
N_{E/\bQ_{17}}  \big(E^\ast\big) = \langle -17 \rangle \cdot \{\ep \in \fo_{17}^\ast \ ; \  \ep \text{ is a $4^{th}$ power modulo $17$} \}.
\]
It is the argument of Section \ref{sec:absenceRPBM} that we will apply below to adelic sections in order to find a Brauer--Manin obstruction for sections of the affine Reichardt--Lind curve.


\subsection{The \texorpdfstring{$\mu_2$}{mu$_2$}-extension}
Let $\pi_1(\cU_{\ov{s}})$ (resp.\ $\pi_1(\cX_{\ov{s}})$) be the fundamental group of the geometric generic fibre of $\cU \to S$ (resp.\ $\cX \to S$) with the base $S=\Spec \bZ[\bruch{1}{34}]$.
The kernel of the quotient map $\pi_1 \surj \pi_1^2$ to the maximal pro-$2$ quotients are characteristic and thus yield short exact sequences
\begin{equation} \label{eq:sesU}
1 \to \pi_1^2(\cU_{\ov{s}}) \to \pi_1^{(2)}(\cU) \to \pi_1(S,\ov{s}) \to 1
\end{equation}
\begin{equation} \label{eq:sesX}
1 \to \pi_1^2(\cX_{\ov{s}}) \to \pi_1^{(2)}(\cX) \to \pi_1(S,\ov{s}) \to 1.
\end{equation}
with the geometrically pro-$2$ fundamental groups $\pi_1^{(2)}(\cU)$ (resp. $\pi_1^{(2)}(\cX)$).
Because the base $S$ and the geometric fibres are $\rK(\pi,1)$ for the prime $2$ and because there is a natural comparison morphism between the Hochschild--Serre spectral sequence for (\ref{eq:sesU}) (resp. (\ref{eq:sesX})) and the Leray spectral sequence for the projection to $S$, we find the following proposition.

\begin{prop}
Let $M$ be a locally constant \'etale sheaf on $\cU$ (resp.\ $\cX$) with fibre a finite abelian $2$-group.
Then the natural maps 
\[
\rH^q(\pi_1^{(2)}(\cU),M) \to \rH^q(\cU,M)
\]
\[
\rH^q(\pi_1^{(2)}(\cX),M) \to \rH^q(\cX,M)
\]
are isomorphisms for all $q\in \bN$. \hfill $\square$
\end{prop}
It follows that the class $\alpha \in \rH^2(\cX,\mu_2)$ of Section \ref{sec:absenceRPBM} already lives in $\rH^2(\pi_1^{(2)}\cX,\mu_2)$ and restricts/inflates to $\alpha_\cU=\chi_{Y} \cup \chi_{17} \in \rH^2(\pi_1^{(2)}\cU,\mu_2^{\otimes 2})$ with $\chi_Y,\chi_{17} \in \rH^1(\pi_1^{(2)}\cU,\mu_2) =  \rH^1(\cU,\mu_2)$ and the identification $\mu_2=\mu_2^{\otimes 2}$.


\subsection{Sections for the Reichardt--Lind curve} \label{sec:sectionsRLcurve}

\begin{thm} \label{thm:nosecRL}
The fundamental group extension $\pi_1(U/\bQ)$ for the affine Reichardt--Lind curve $U/\bQ$ does not split. In particular, the section conjecture holds trivially for $U/\bQ$ as there are neither rational points nor sections.

More precisely, the geometrically pro-$2$ extension $\pi_1^{(2)}(X/\bQ)$ of the projective Reichardt--Lind curve $X/\bQ$ does not admit a section $s$ that allows  locally at $p=2$ and $p=17$ a lifting $\tilde{s}_p$
\[
\xymatrix@M+1ex@R-2ex{& \Gal_{\bQ_p} \ar@{.>}[dl]_{\tilde{s}_p} \ar[d]^{s \otimes \bQ_p} \\
\pi_1^{(2)}(U) \ar@{->>}[r] & \pi_1^{(2)}(X) .}
\]
\end{thm}

Note that the method of \cite{hs} does not apply to the Reichardt--Lind curve because its jacobian $E$ has rank $0$, see the end of Section \ref{sec:geoarithRL}.

\begin{pro}
We argue by contradiction. Let $s$ be a section of $\pi_1^{(2)}(X/\bQ)$ that allows local lifts $\tilde{s}_p$ at $p=2$ and $p=17$. We will compute the Brauer--Manin obstruction $\langle \alpha, (s \otimes \bQ_v)\rangle$. 

The section  $s$ descends to a section $\sigma : \pi_1(S) \to \pi_1^{(2)}(\cX)$ with $S=\Spec \bZ[\bruch{1}{34}]$ as above by Section \ref{sec:notameram}. It follows that  $(s \otimes \bQ_v)^\ast \alpha$ is the image of $\sigma^\ast(\alpha)$ under the map 
\[
\rH^2(\pi_1 S,\mu_2) \subset \rH^2(S,\mu_2)  \to \rH^2(\bQ_v,\mu_2) = \ _2\Br(\bQ_v)
\]
and thus vanishes for $v \not= 2,17,\infty$. It remains to compute the contribution to $\langle \alpha, (s \otimes \bQ_v)_v \rangle$ at the places $v=2,17$ and $\infty$. 

When restriting to $\bR=\bQ_v$ for $v=\infty$ we know from the real section conjecture, see \cite{stix:periodindex} Appendix A, that $s \otimes \bR$ belongs to a real point $x \in X(\bR)$, that we moreover can move freely in its connected component, so that it is in $U(\bR)$. The section $s \otimes \bR$ thus lifts to a  $\tilde{s}_\infty : \Gal_\bR \to \pi_1^{(2)}(U)$.
We conclude that  for $v=2,17$ and $\infty$
\[
(s \otimes \bQ_v)^\ast(\alpha) = \tilde{s}_v^\ast(\alpha_U) = \tilde{s}_v^\ast(\chi_Y \cup \chi_{17}) = \tilde{s}_v^\ast(\chi_Y) \cup \tilde{s}_v^\ast(\chi_{17}) = \chi_{Y(\tilde{s}_v)} \cup \chi_{17}
\]
which is the class of the quaternion algebra  $(y_v,17) \in \Br(\bQ_v)$ with $y_v=Y(\tilde{s}_v) \in \bQ_v^\ast/(\bQ_v^\ast)^2$. Because $17$ is a square at $2$ and $\infty$, the local contribution at $v=2$ and  at $v=\infty$ vanishes.

As $X^4-17$ is a norm for the projection $U \times_\bQ \bQ(\sqrt[4]{17}) \to U$ and equal to the unit $2Y^2$ we find by Section \ref{sec:evNorm} that $2y_{17}^2$ is a norm from $\bQ_{17}(\sqrt[4]{17})$. At first the equation
\[
2y_{17}^2 = N_{ \bQ_{17}(\sqrt[4]{17})/\bQ_{17}}\Big(\big(X-\sqrt[4]{17}\big)(\tilde{s}_{17})\Big)
\]
holds only in $\bQ_{17}^\ast/(\bQ_{17}^\ast)^2$. But with a suitable choice of representative $y$ for $y_{17} \in \bQ_{17}^\ast/(\bQ_{17}^\ast)^2$  we find that $2y^2$ is a norm from $\bQ_{17}(\sqrt[4]{17})$ already in $\bQ_{17}^\ast$. We conclude as at the end of Section  \ref{sec:absenceRPBM} that $y$ is not a norm from $\bQ_{17}(\sqrt{17})$ and thus $\inv_{17}\big((y,17)\big) = \bruch{1}{2}$ is nontrivial, so that also 
\[
\langle \alpha, (s \otimes \bQ_v)\rangle = \sum_v \inv_v (s\otimes \bQ_v)^\ast(\alpha) = \inv_{17}\big((y,17)\big) = \bruch{1}{2} \not= 0
\]
which contradicts Theorem \ref{thm:globalsectionsBMkernel}.
\end{pro}


\subsection{Arithmetic twists of the Reichardt--Lind curve}
In this section we study an arithmetically twisted version of the Reichardt--Lind curve. Let $k$ be an algebraic number field and let $X$ be the smooth projective curve over $k$ defined by the homogeneous equations 
\begin{equation} \label{eq:twistedRL}
\left\{\begin{array}{l}
\ell T^2 = A^2 - p B^2 \\
AB = C^2
\end{array}\right.
\end{equation}
with $\ell,p \in k^\ast$. The curve $X$ has good reduction for places $v \nmid 2\ell p$.  Let $U \subset X$ be the complement  of the support $D = P \amalg Q$ of $\divisor(T/B)$, with $P \cong \Spec(k[t]/t^4-p)$ and $Q= \Spec(k[t]/t^2-\ell)$. The Brauer class $\alpha_U = \chi_{T/B} \cup \chi_p \in \rH^2(U,\mu_2)$  lifts to a class $\alpha \in \rH^2(X,\mu_2)$ because its residues at $P$ and $Q$  vanish as $T/B$  (resp.\ $p$) becomes a square near $Q$ (resp.\ $P$).

For the computation of the Brauer--Manin obstruction with respect to $\alpha$ along the lines of  Section \ref{sec:sectionsRLcurve}  we impose the following list of conditions.
\begin{itemize}
\item[(i)] $p$ is positive with respect to any real place of $k$.
\item[(ii)] for all $v \mid 2\ell p$ with $\sqrt{p} \not \in k_v$ we have $\mu_4 \subset k_v$ and
\[
N=\# \{ v \mid 2 \ell p \ ; \ \sqrt{p} \not \in k_v \text{ and $\ell$ not a norm from } k_v(\sqrt[4]{p})/k_v \}
\]
is odd.
\end{itemize}
We furthermore need lifts $\tilde{s}_v:\Gal_{k_v} \to \pi_1^{(2)}(U)$ locally at places $v \mid 2 \ell p$, but then  the conditions (i) \& (ii) enforce $\langle \alpha, (s_v)\rangle = N/2= 1/2$ regardless of the choice of the adelic section $(s_v)$. 
The curve $X$ has local points for all places $v$ of $k$ if 
\begin{itemize}
\item[(iii)] for all $v \mid 2\ell p$  we have $\sqrt{\ell} \in k_v$ or $\sqrt[4]{p} \in k_v$.
\end{itemize}
because then $P$ or $Q$ split locally for $v \mid 2 \ell p$. The places of good reduction (resp.\ the infinite places) have local points anyway (resp.\ by (i)).

In order to get a simpler list of conditions we specialize now to the case where $p$ is a prime element of $\fO_k$ that is coprime with $2\ell$. Moreover, we choose in (ii) for all places $v \not= p$ that $\sqrt{p} \in k_v$. It follows that the local quadratic norm residue symbol $(\ell,p)_v$ vanishes for all places $v$ but possibly $v=p$, hence also at $v=p$ by reciprocity and thus $\sqrt{\ell} \in k_v$ for $v=p$. Using local class field theory to determine the norm group of $k_v(\sqrt[4]{p})/k_v$ for $v=p$, it is thus enough to ask
\begin{itemize}
\item[(i')] $p$ is a prime element coprime to $2 \ell$ and positive with respect to any real place of $k$.
\item[(ii')] for $v=p$ we have $\mu_4 \subset k_v$ and $\ell$ is not congruent to a $4^{th}$ power modulo $p$.
\item[(iii')] for all $v \mid 2\ell$  we have $\sqrt{p} \in k_v$ and a local solution in $k_v$, e.g., by requiring  $\sqrt[4]{p} \in k_v$.
\end{itemize}
We further simplify to $k=\bQ$ and set $\ell=\ep \prod_{i=1}^n q_i$ as a product of $n\geq1$ positive distinct prime numbers $q_i$ and a sign $\ep = \pm 1$. Moreover we achieve condition (iii')  for $v \nmid 2$ by asking $\sqrt{p} \in k_v$ but $\mu_4 \not\subset k_v$. Hence we need to satisfy the conditions
\begin{itemize}
\item[(i'')] $p>0$ is an odd prime number different from the prime factors $q_i$ of $\ell$.
\item[(ii'')] $p \equiv 1 \mod 4$ and $\ell$ is not congruent to a $4^{th}$ power modulo $p$.
\item[(iii'')]  $p$ is a square modulo $q_i$ for all $i =1, \ldots,n$, and $q_i \equiv 3 \mod 4$ for $q_i$ odd. 
\item[(iv'')] $p \equiv 1 \mod 8$ and we have a local solution in $\bQ_2$, e.g., by requiring  $p \equiv 1 \mod 16$.
\end{itemize}
Let $\ell$ be even (but not divisible by $4$). The number field $F=\bQ(\sqrt[4]{\ell},\sqrt{q_1},\ldots,\sqrt{q_n},\mu_8)$ is of degree $2$ over $E=\bQ(\sqrt{q_1},\ldots,\sqrt{q_n},\mu_8)$. Conditions (i'') - (iii'')  on $p$ translates to the condition that $p$ is unramified in $F/\bQ$ and the Frobenius $\Frob_\fp$ of a place $\fp| p$ for the extension $F/\bQ$ generates the central group $\Gal(F/E)$. The solution in $\bQ_2$ is now automatic as follows. We write the affine version of equation (\ref{eq:twistedRL}) as $Z^4 = p + \ell Y^2$, and $p+\ell Y^2$ must be $\equiv 1 \mod 16$ in order to allow a $4^{th}$ root in $\bQ_2$.
If $p \equiv 1 \mod 16$ we put $Y=0$, if $p \equiv 9 \mod 16$ we put $Y=2$, and we are done. Using the Chebotarev density theorem and $\sqrt{2} = (1+i)/\zeta_8$ we find that a set of density $\bruch{1}{2^{n+2}}$ of all prime numbers yield valid choices for $p$ when $\ell$ is even and squarefree and all odd prime factors of $\ell$ are $\equiv 3 \mod 4$.

In case $\ell$ is odd, we have to replace $F/E$ by $\bQ(\sqrt[4]{\ell},\sqrt{q_1},\ldots,\sqrt{q_n},\mu_{16})/\bQ(\sqrt{q_1},\ldots,\sqrt{q_n},\mu_{16})$ with the same translation of the conditions on $p$ 
to the Frobenius $\Frob_\fp$ for a place $\fp | p$. In this case, the Chebotarev density theorem yields a set of density $\bruch{1}{2^{n+4}}$ of all prime numbers as valid choices for $p$ when $\ell$ is odd and squarefree and all prime factors of $\ell$ are $\equiv 3 \mod 4$.
Concrete examples are given as follows.
\begin{enumerate}
\item $\ell=2$ and $p=17$, the Reichardt--Lind curve, or $p=41$ or $p=97$,
\item $\ell=6$ and $p=73$,
\item $\ell=11$ and $p=97$,
\item $\ell=19$ and $p=17$.
\end{enumerate}
All these arithmetic twists yield  more empty examples for the section conjecture where the existence of sections is Brauer--Manin obstructed. In particular, there are infinitely many of these examples.


\subsection{An isotrivial family of examples} In the introduction of \cite{poonen:familyHasseP} Poonen reports an example constructed by Elkies. The smooth fibres of the isotrivial family $X_t$ with $t \in \bP^1$ of arithmetic twists of the Reichardt-Lind curve, given in affine form by the equation
\[
2Y^2 = Z^4 - N(t), \qquad
N(t) = \left(1 + \bruch{2}{1+t+t^2}\right)^4 + 16,
\]
yield for each $t \in \bQ \cup \{\infty\}$ a counterexample to the Hasse principle as explained below. Note that for $t = \infty$ the fibre $X_\infty$ coincides with the original Reichardt--Lind curve.

\subsubsection{Local points}
The key property of $N(t)$ is that it takes only odd values at rational $t$. After rescaling we may replace $N=N(t)$ by a quartic free integer $N_0 = M^4 \cdot N(t)$ which equals $N_0=A^4 + 16B^4$ with $A,B$ odd and coprime integers.  Clearly we have $\bR$-points as $N_0>0$ and also $\bQ_p$-points for $p \nmid 2N_0$, the primes of good reduction of $X_t$. 

For local points at $p=2$ we note that $N_0 \equiv 1 \mod 16$, and thus $N_0$ is a 4$^{th}$ power in $\bQ_2$, which leads to a $\bQ_2$-point with $Y=0$. For an odd prime $p \mid N_0$ we find that $-1$ is a 4$^{th}$ power in $\bF_p$, hence $p \equiv 1 \mod 8$, and so $2$ is a square in $\bF_p$. It follows that there is $y \in \bZ_p$ such that $N_0+2y^2$ is a 4$^{th}$ power modulo $p$ and thus also a 4$^{th}$ power in $\bQ_p$ which yields a $\bQ_p$-point.

\subsubsection{The absence of sections}
Let $U_t$ be the open subscheme of $X_t$ given by 
\[
2Y^2 = Z^4 - N(t), \qquad Y \not= 0. 
\]
\begin{prop}
The  fundamental group extension $\pi_1(U_t/\bQ)$ for the arithmetic twist $U_t$ with $t \in \bQ$ of the affine Reichardt--Lind curve does not allow a section.
\end{prop}
\begin{pro}
We rescale as above and replace $N$ by $N_0 = A^4+16B^4$ with $A,B$ odd and coprime.
The Brauer--Manin obstruction  for an adelic section $(s_p)$ with respect to the Brauer class $\alpha = \chi_Y \cup \chi_{N_0}$ can have a nontrivial local contribution 
\[
\inv_p(s_p^\ast(\alpha)) = \inv_p\big((Y(s_p),N_0)\big)
\]
 at most at places $p \mid N_0$. Then $p \equiv 1 \mod 8$ and $\bQ_p$ contains $\mu_4$. Hence $\bQ_p(\sqrt[4]{N_0})$ either equals $\bQ_p$ and then $\inv_p(s_p^\ast(\alpha)) = 0$, or $\bQ_p(\sqrt[4]{N_0})$ is an abelian cyclic extension of $\bQ_p$ of degree $4$. Local class field theory shows, as in Section \ref{sec:absenceRPBM}, that multiplication by $2$ induces an injective map
\[
\mu_2 \cong  \bQ_p^\ast/ \rN_{\bQ_p(\sqrt{N_0})/\bQ_p} \big(\bQ_p(\sqrt{N_0})^\ast\big)  \ \inj \ \bQ_p^\ast/ \rN_{\bQ_p(\sqrt[4]{N_0})/\bQ_p} \big(\bQ_p(\sqrt[4]{N_0})^\ast\big)  \cong \mu_4.
\]
As $2Y(s_p)^2$ is a norm from $\bQ_p(\sqrt[4]{N_0})$, we conclude that $\inv_p(s_p^\ast(\alpha))$ is nontrivial, namely $Y(s_p)$ is not a norm form $\bQ_p(\sqrt{N_0})$, if and only if $2$ is not a norm from 
$\bQ_p(\sqrt[4]{N_0})$. If  the exponent $v_p(N_0)$ of $p$ in $N_0$ is even, then $2$ is not a norm from $\bQ_p(\sqrt[4]{N_0})$ if and only if $2$ is not a square in $\bF_p$, but which it always is due to $p \equiv 1 \mod 8$.

If  the exponent $v_p(N_0)$ of $p$ in $N_0$ is odd, then $2$ is not a norm from $\bQ_p(\sqrt[4]{N_0})$ if and only if $2$ is not a 4$^{th}$ power in $\bF_p$, hence for such primes $p$ we find
\[
\inv_p(s_p^\ast(\alpha)) = \left(\bruch{2}{p}\right)_4
\]
under the identification $\bruch{1}{2}\bZ/\bZ = \mu_2(\bF_p)$, and where $\left(\bruch{}{p}\right)_4$ is the quartic residue map
\[
\left(\bruch{}{p}\right)_4 \ : \ \bF_p^\ast \to \mu_4(\bF_p), \qquad x \mapsto \left(\bruch{x}{p}\right)_4 = x^{(p-1)/4}.
\]
The Brauer--Manin pairing thus evaluates as
\[
\langle \alpha, (s_p)\rangle = \#\{ p \mid N_0 \ ; \ v_p(N_0) \text{ odd, and } \left(\bruch{2}{p}\right)_4 = -1 \} \cdot \bruch{1}{2} \in \bQ/\bZ.
\]
Let $p$ be a prime dividing $N_0$, hence of the form $p=a^2 + 16b^2$ with $a$ odd. Let $\zeta_p$ be a $p^{th}$ root of unity. The subfield $k$ of $\bQ(\zeta_p)$ of degree $4$ over $\bQ$ is uniquely determined by being cyclic of degree $4$ over $\bQ$ with ramification only above $p$.  Using \cite{jly} Thm 2.2.5, we can determine $k = \bQ(\lambda)$ with $\lambda^2 = \sqrt{p}(4b + \sqrt{p}).$
The decomposition behaviour above $(2)$ is given on the one hand by the Frobenius element 
\[
\left(\bruch{2}{p}\right)_4 \in \mu_4(\bF_p) = \Gal(k/\bQ),
\]
and on the other hand explicitly as follows. The prime $(2)$ is completely split in $k$ if and only if $\sqrt{p}(4b + \sqrt{p})$ is a square in $\bQ_2$, or equivalently 
\[
\sqrt{p}(4b + \sqrt{p}) \equiv 1 \mod 8
\]
with respect to a choice of $\sqrt{p} \in \bQ_2$. As $p \equiv 1  \mod 8$, the latter is equivalent to $b$ being even.

We resort to arithmetic of $\bZ[i]$ to prove that any way of writing $N_0$ as $(A^2)^2 + 16(B^2)^2$ comes from writing its prime factors, which are all congruent to $1$ modulo $8$, in the form $p=a^2+16b^2$ for suitable $a,b \in \bZ$, and then making inductively use of the identity
\begin{equation} \label{eq:normidentity}
(a^2 +16b^2)\cdot (c^2 + 16 d^2) = \big(ac-16bd\big)^2 + 16 \big(ad+bc\big)^2.
\end{equation}
In our application of (\ref{eq:normidentity}) the integers  $a,c$ are always odd. Thus the parity of the second component $ad+bc \equiv b+d \mod 2$ is additive. We conclude that modulo $2$
\[
\#\{ p \mid N_0 \ ; \ v_p(N_0) \text{ odd, and } \left(\bruch{2}{p}\right)_4 = -1 \} 
\ \equiv \sum_{p \mid N_0 \ ;   \left(\bruch{2}{p}\right)_4 = -1} v_p(N_0)  \ \equiv \ B^2 \ \equiv \ 1 ,
\]
and thus $\langle \alpha, (s_p)\rangle \not= 0$. Hence $\pi_1(U_t/\bQ)$ does not admit sections due to the Brauer--Manin obstruction coming from $\alpha$.
\end{pro}


\section{Beyond sections of the abelianized fundamental group}
The known obstructions against sections use for example period and index in the $p$-adic local case or the Brauer--Manin obstruction  in Section \ref{sec:BM} and Section \ref{sec:RL}. 
The Brauer--Manin obstruction extends to $0$-cycles. 
As an obstruction against $0$-cycles of degree $1$ it 
only depends on abelian information that takes place on the generalized semiabelian Albanese torsor $U \to \SAlb_{U/k}^1$ of the curve $U/k$. On curves of genus $\geq 2$, however, the Brauer--Manin obstruction depends on information beyond the abelian level, as there are examples of curves of index $1$ with adelic points but no global point, see \cite{sko:torsors} p.128 on an example of \cite{bremnerlewismorton}.

Nevertheless, as the Reichardt--Lind curve has genus $1$, the natural question arises, whether the above method actually prove that the geometrically abelianized fundamental group extension $\pi_1^\ab(U/k) = \pi_1(\SAlb_{U/k}^1/k)$ denies having sections under the circumstances that prove the nonexistence of sections for the full extension $\pi_1(U/k)$. We will discuss this issue first for the above proof in the case of the Reichardt--Lind curve and then give an example for the local $p$-adic obstruction on a curve of genus $2$.


\subsection{Abelian is enough for the Reichardt--Lind curve}
We keep the notation from the preceding Section \ref{sec:RL}. Let $\pi_1^{(\ab,2)}(U/\bQ)$ be the maximal geometrically abelian pro-$2$ quotient of $\pi_1(U/\bQ)$. 
\begin{lem} \label{lem:abel}
(1) Inflation induces an isomorphism $\rH^1(\pi_1^{(\ab,2)}(U),\mu_2) \to \rH^1(\pi_1^{(2)}(U),\mu_2)$.

(2) The class $\alpha_U = \chi_Y \cup \chi_{17}$ lifts uniquely to a class $\alpha_U^\ab= \chi_Y \cup \chi_{17}$, where $\chi_Y$ and $\chi_{17}$ are the corresponding classes under the isomorphism from (1).
\end{lem}
\begin{pro}
The Hochschild--Serre spectral sequence $\ ^\ab\rE$ for $\pi_1^{(\ab,2)}(U) \to \Gal_{\bQ}$ maps via inflation to the corresponding spectral sequence $\rE$ for $\pi_1^{(2)}(U) \to \Gal_{\bQ}$. The map $\ ^\ab\rE_2^{p,q} \to \rE_2^{p,q}$ is an isomorphism for $q=0,1$ which implies (1) and thus (2).
\end{pro}

For the argument of the proof of Theorem \ref{thm:nosecRL} to work with local lifts to the abelianized fundamental group extension $\pi_1^{(\ab,2)}(U \times_\bQ \bQ_v/\bQ_v)$ for $v=2,17$ we would need that the lift $\alpha$ of $\alpha_U$ maps to $\alpha_U^\ab$ under 
\[
\rH^2(\pi_1^{(2)}(X),\mu_2) \to \rH^2(\pi_1^{(\ab,2)}(U),\mu_2) \to \rH^2(\pi_1^{(2)}(U),\mu_2).
\]
Let $\rF^\bullet\rH^\ast$ denote the filtration on the cohomology groups coming from the Hochschild--Serre spectral sequences with respect to the projection to $\Gal_\bQ$. The proof of Lemma \ref{lem:abel} shows that $\rF^1\rH^2(\pi_1^{(\ab,2)}(U),\mu_2) \to \rF^1\rH^2(\pi_1^{(2)}(U),\mu_2)$ is an isomorphism. The classes $\alpha_U$ and $\alpha_U^\ab$ die under the projection to $\gr_{\rF}^0\rH^2$ due to the constant factor $\chi_{17}$ and thus lie in the respective $\rF^1\rH^2$.
\begin{lem} \label{lem:gr0inj}
The natural map $\gr^0_{\rF} \rH^2(\pi_1^{(2)}(X),\mu_2) \to \gr^0_{\rF} \rH^2(\pi_1^{(\ab,2)}(U),\mu_2)$ is injective.
\end{lem}
\begin{pro}
The map injects into the map $\rH^2(\pi_1^{2}(X\times_\bQ \bQ^\alg),\mu_2) \to \rH^2(\pi_1^{\ab,2}(U\times_\bQ \bQ^\alg),\mu_2)$ which is split injective because the surjection 
\[
\pi_1^{\ab,2}(U\times_\bQ \bQ^\alg) \surj \pi_1^{2}(X\times_\bQ \bQ^\alg)
\]
admits a section.
\end{pro}

As a consequence of Lemma \ref{lem:gr0inj} we deduce that a lift $\alpha$ that is suitable for a version with only abelian local lifts $\tilde{s}_v$ at $v=2,17$ has to lie in $\rF^1\rH^2(\pi_1^{(2)}(X),\mu_2)$. 

\begin{prop} \label{prop:goodlift}
The class $\alpha_U^\ab \in  \rH^2(\pi_1^{(\ab,2)}(U),\mu_2)$ admits a lift $\alpha$ to $\rH^2(\pi_1^{(2)}(X),\mu_2)$.
\end{prop}
\begin{pro}
Let $\ov{X}$ (resp.\ $\ov{U}$) denote the base change to the algebraic closure of $X$ (resp.\ $U$).
As $\rF^2\rH^2(\pi_1^{(2)}(X),\mu_2)$ equals $\rF^2\rH^2(\pi_1^{(\ab,2)}(U),\mu_2)$ we have to study the map 
\[
\gr^1_{\rF} \rH^2(\pi_1^{(2)}(X),\mu_2) \to \gr^1_{\rF} \rH^2(\pi_1^{(\ab,2)}(U),\mu_2)
\]
which agrees with 
\[
\rH^1\big(\bQ,\rH^1(\ov{X},\mu_2)\big) \to \rH^1\big(\bQ,\rH^1(\ov{U},\mu_2)\big)
\]
by the following lemma.
\begin{lem}
The maps $\rE_2^{3,0} \to \rE_\infty^3$ for the Hochschild--Serre spectral sequence for the maps  $\pi_1^{(2)}(X) \to \Gal_\bQ$ (resp.\ of $\pi_1^{(\ab,2)}(U/\bQ) \to \Gal_Q$) are injective.
\end{lem}
\begin{pro}
By Tate's theorem we have $\rE_2^{3,0} = \rH^3(\bQ,\mu_2) = \rH^3(\bR,\mu_2)$, so that we can compute the injectivity of the map $\rE_2^{3,0} \to \rE_\infty^3$ by first base changing to $\bR$. There our curve has a point, which yields a retraction of the edge map showing it to be injective.
\end{pro}

We continue with the proof of Proposition \ref{prop:goodlift}.
Recall that $D=X-U$ is the boundary divisor, and that canonically $\mu_2^{\otimes 2} = \mu_2$.
The image of $\alpha_U^\ab$ in $\rH^1\big(\bQ,\rH^1(\ov{U},\mu_2)\big)$ is the image of $\chi_{17}$ under $\rH^1(\bQ,\mu_2 \otimes -)$ applied to the map $\bZ/2\bZ \to \rH^1(\ov{U},\mu_2)$ given by $1 \mapsto \chi_Y$. For deciding the existence of a preimage we use the exact sequence
\[
0\to \rH^1(\ov{X},\mu_2^{\otimes 2}) \to \rH^1(\ov{U},\mu_2^{\otimes 2}) \xrightarrow{\res} \mu_2[D(\bQ^\alg)]^{\sum\equiv 0} \to 0
\]
that maps $\res(\chi_Y)$ to the sum of the conjugates of $P-2Q$, hence modulo $2$ just the conjugates of $P$. In the resulting cohomology sequence 
\[
\rH^1(\bQ,\rH^1(\ov{X},\mu_2^{\otimes 2})) \to \rH^1(\bQ,\rH^1(\ov{U},\mu_2^{\otimes 2}) ) \xrightarrow{\res} \rH^1(\bQ,\mu_2[D(\bQ^\alg)]^{\sum\equiv 0})
\]
the image of $\alpha_U^\ab$ in the right group is the image of $\chi_{17}$ under the map 
\[
\rH^1(\bQ,\mu_2) \to \rH^1(\bQ,\mu_2[D(\bQ^\alg)]^{\sum\equiv 0})
\]
induced by $1 \mapsto \res(\chi_Y)$. The short exact sequence
\[
0 \to \mu_2[D(\bQ^\alg)]^{\sum\equiv 0} \to \mu_2[D(\bQ^\alg)] \xrightarrow{\sum} \mu_2 \to 0
\]
yields a cohomology sequence
\[
0 \to \mu_2 \xrightarrow{\delta} \rH^1(\bQ,\mu_2[D(\bQ^\alg)]^{\sum\equiv 0}) \to \bQ(\sqrt[4]{17})^\ast/\big(\bQ(\sqrt[4]{17})^\ast\big)^2 \times  \bQ(\sqrt{2})^\ast/\big(\bQ(\sqrt{2})^\ast\big)^2.
\]
The image of $17 =\chi_{17} \in \bQ^\ast/(\bQ^\ast)^2=\rH^1(\bQ,\mu_2)$ in the right group by the map induced from $1 \mapsto \res(\chi_Y)$ vanishes, so that the image of $\alpha_U^\ab$ lies in the image of $\delta$. Restricting to $\bQ(\sqrt{17})$ decomposes $D$ into three divisors of degree $2$, hence the map $\delta$ remains injective, although $\chi_{17}$ dies, the conclusion of which is, that the image of $\chi_{17}$ and thus $\alpha_U^\ab$ in 
$\rH^1(\bQ,\mu_2[D(\bQ^\alg)]^{\sum\equiv 0})$ must vanish. 
\end{pro}

\begin{cor} \label{cor:abelianRL}
The abelianized fundamental group extension $\pi_1^\ab(U/\bQ)$ for the affine Reichardt--Lind curve $U/\bQ$ does not split.

More precisely, the geometrically pro-$2$ extension $\pi_1^{(2)}(X/\bQ)$ of the projective Reichardt--Lind curve $X/\bQ$ does not admit a section $s$ that allows  a lifting $\tilde{s}_p$
\[
\xymatrix@M+1ex@R-2ex{& \Gal_{\bQ_p} \ar@{.>}[dl]_{\tilde{s}_p} \ar[d]^{s \otimes \bQ_p} \\
\pi_1^{(\ab,2)}(U) \ar@{->>}[r] & \pi_1^{(2)}(X) }
\]
locally at $p=2$ and $p=17$. 
\end{cor}

In view of  a recent result of Esnault and Wittenberg \cite{ew:birationalabelian}, that an abelian birational section implies a divisor of degree $1$, our Corollary \ref{cor:abelianRL} is merely an explicit form showing `how birational' the section must be assumed to be. Indeed, the projective  Reichardt--Lind curve is a nontrivial torsor under its jacobian and thus has index different from $1$. Therefore $X/\bQ$ does not allow a birational abelian section by \cite{ew:birationalabelian} Theorem 2.1. Note that the latter is conditional on the finiteness of $\Sha^1(\bQ,\Pic_X^0)$ which is known for the Reichardt--Lind curve.


We close the discussion of the Reichardt--Lind curve with the following natural question.

\begin{ques}
Does the extension $\pi_1(X/\bQ)$ for the projective Reichardt--Lind curve $X$ split?
\end{ques}


\subsection{An example in genus \texorpdfstring{$2$}{2}} The idea is to exploit the annoying factor $2$ from Lichtenbaum's congruences for period and index of a smooth projective curve over a $p$-adic field which shows up in the main result of \cite{stix:periodindex} Theorem 16.

\begin{thm}
Let $k$ be a finite extension of $\bQ_p$, and let $X/k$ be a smooth projective curve of genus at least $1$. 
\begin{itemize}
\item[(1)] If $\pi_1^\ab(X/k)$ splits, then $\pe(X)$ is a power of $p$. 
\item[(2)] As a partial converse: if $\pe(X)$ equals $1$, then $\pi_1^\ab(X/k)$ splits.
\item[(2)]  If we have a section for the natural map to $\Gal_k$ of the maximal quotient of $\pi_1(X)$ which geometrically  is the product of an abelian group of order prime to $2$ with a $2$-group, which is an abelian extension of the $\rH_1(X \times_k k^\alg, \bZ/2\bZ)$-quotient, then the index is also a power of $p$.
\end{itemize}
\end{thm}
\begin{pro}
This follows from a more careful study of \cite{stix:periodindex} Section 5.
\end{pro}

In order to get a curve $X/k$ such that $\pi_1^\ab(X/k)$ splits despite the fact that $\pi_1(X/k)$ does not allow sections, we are led to look for curves of genus $2$, period $1$ and index $2$. These are guaranteed by results of Clark \cite{clark} and Sharif \cite{sharif}. 

Here is an alternative construction for $p \not= 3$. Let $P$ be the Brauer--Severi variety for the Azumaya algebra of invariant $1/2$ over $k$. The curve $P$ has genus $0$ and thus has a point $y$ of degree $2$. There is  a non-constant rational function $f$ on $P$ with divisor $z-y$, where $z$ is then necessarily also a point of degree $2$. The branched cover $h:X \to P$ defined by a cubic root of $f$ is  totally branched above $y,z$, and unramified elsewhere. The curve $X$ is of genus $2$ by the Riemann--Hurwitz formula, its index is different from $1$ but divides $2g-2$, hence $\ix(X) = 2$, and the period divides $g-1$, hence equals $1$. The curve $X/k$ therefore gives an explicit example for a curve such that $\pi_1^\ab(X/k)$ splits although $\pi_1(X/k)$ does not. 

\begin{ques}
Find an explicit example for a smooth projective curve $X$ over a number field $k$ which is a counter example to the Hasse principle and such that $\pi_1^\ab(X/k)$ splits despite the fact that $\pi_1(X/k)$ does not allow sections. The existence of adelic points rules out that $\pi_1(X/k)$ does not split for local reasons and thus may be weakened to the requirement that there are local sections everywhere for still an interesting example.
\end{ques}


\section{Example: the Selmer curve} \label{sec:selmer}


\subsection{Geometry and arithmetic of the Selmer curve}  
Another famous example of a curve violating the Hasse principle was found by Selmer in \cite{selmer:curve} as the plane cubic curve $S/\bQ$ of genus $1$ which in homogeneous coordinates is given by
\[
3X^3 + 4Y^3 + 5Z^3 = 0.
\]
The Jacobian $E=\Pic_S^0$ of $S$ is an elliptic curve with complex multiplication by $\bZ[\zeta_3]$ given by the homogeneous equation
\[
A^3 + B^3 + 60C^3 = 0,
\]
and $[1:-1:0]$ as its origin. The curve $S$, as a principal homogeneous space under $E$, describes a nontrivial $3$-torsion element of $\Sha^1(\bQ,E)$, see \cite{mazur:localtoglobal} I \S4+9. 

Cassels gives in \cite{cassels:lectures} a proof of the absence of global points on $S$ via \'etale covers as follows. There is a finite \'etale cover $h : S \to E$ of degree $9$, which is a torsor under the finite \'etale group scheme $E[3]$. Standard methods in the arithmetic of elliptic curves allow Cassels to show that $E(\bQ) = 0$, see \cite{cassels:lectures}  \S18 Lem 2, so that $S(\bQ)$ is contained in $h^{-1}(0) = S \cap \{XYZ = 0\}$, which does not  contain rational points by direct inspection. A sketch of Selmer's original proof based on a $3$-descent can be found in \cite{cassels:diopheq} p.\ 206.

As $\Sha^1(\bQ,E)$ is finite, see \cite{mazur:localtoglobal} Thm 1 and \S9, there must be a locally constant Brauer--Manin obstruction on $S$ which explains the absence of rational points. We proceed with constructing the appropriate global Brauer class in the next sections. The construction relies on computations with a computer algebra package (SAGE). As the results on the Selmer curve only illustrate the difficulties with the Brauer--Manin obstruction in the section conjecture, we can afford avoiding the details of these computations.

\subsection{A degree $3$ cover} For the rest of Section \ref{sec:selmer} we set $k=\bQ(\zeta_3)$ for a fixed cubic root of unity $
\zeta_3$, and $K=k(\ep)$ with $\ep = \sqrt[3]{6}$. 
We denote by $\sigma$ the generator of $\Gal(K/k)$ which maps $\ep$ to $\zeta_3 \ep$. The element 
\[
\gamma = \zeta_3 \ep^2 + (2\zeta_3 +1)\ep + 2\zeta_3 \in K
\]
is constructed as an element satisfying $N_{K/k}(\gamma) = -10$ using SAGE, and turns out to be a product of the prime above $2$ by a prime above $5$. We find
\[
N_{K/k}\left(\bruch{2Y+\ep X}{\gamma Z}\right) = 1,
\]
and the Ansatz of the Lagrange resolvent  yields 
\[
U = \bruch{(2Y+\ep X)(2Y+\zeta_3 \ep X) + \gamma Z (2Y+\zeta_3 \ep X) + \gamma \cdot \sigma(\gamma) Z^2}{(2Y+\ep X)(2Y+\zeta_3 \ep X) } \in K(S)^\ast,
\]
such that 
\[
\sigma(U)/U = \bruch{2Y+\ep X}{\gamma Z}
\]
We define $F \in k(S)^\ast$ as
\begin{equation} \label{eq:defF}
F = N_{K/k}(U) = \bruch{\sigma^2(\gamma)}{\gamma} \cdot U^3 \cdot \bruch{2Y + \ep X}{2Y + \zeta^2 \ep X},
\end{equation}
so that its divisor equals
\[
\divisor(F) = 3 \divisor(U) + 3 \cdot ( [-2:\ep:0] -[2,\zeta_3^2\ep:0]) = 3 D.
\]
The divisor $D$ is $k$-rational, of degree $0$ and describes a nontrivial $3$-torsion class in $E(k)$.
\begin{lem} \label{lem:E3}
The $3$-torsion $E[3]$ of $E$ sits in a short exact sequence of $\Gal_\bQ$-modules 
\begin{equation} \label{eq:E3}
0 \to \mu_3 \to E[3] \to \bZ/3\bZ \to 0,
\end{equation}
which as an extension is given by the Kummer character $\chi_{60} \in \rH^1(\bQ,\mu_3)$.
In particular, for any field extension $k'/\bQ$, for which $60$ is not a cube in $k'$, the restriction of (\ref{eq:E3}) to $\Gal_{k'}$-modules does not split and $E[3](k') = \mu_3(k')$.
\end{lem}
\begin{pro}
The group $E[3]$ is given by the solutions of $ABC=0$, and the point $[-1,\zeta_3,0]$ generates a subgroup $\mu_3 \subset E[3]$. The quotient $E[3]/\mu_3$ has to be constant by $\det E[3] \cong \mu_3$ and the rest follows from computing the Galois action on the preimage of a generator of $\bZ/3\bZ$.
\end{pro}

We deduce from Lemma \ref{lem:E3}, that $D$ generates $\mu_3 \subset E[3]$. The associated $3$ isogeny
\[
f : E' \to E
\]
has Galois group $\bZ/3\bZ$ , the Cartier dual of the group generated by $D$. The cover 
\[
h: W \to S_k = S \times_\bQ k
\]
which is constructed as the integral closure in the extension defined by $\sqrt[3]{F}$ yields a twisted version of $f:E' \to E$. The associated surjective group homomorphism
\[
m_F : \pi_1(S_k) \to \bZ/3\bZ \quad \in \rH^1(\pi_1(S_k), \bZ/3\bZ)
\]
computes as $\zeta_3^{m_F(g)} = \chi_F(g)$ for $g \in \pi_1(S_k)$ and where $\chi_F$ is the cubic Kummer character associated to $F$. The restriction to the pro-$3$ component of the geometric fundamental group we denote by 
\[
\ph = m_F|_{\pi_1(\bar{S})} : \rT_3 E = \pi_1^{3} (\ov{E}) = \pi_1^3 (\ov{S}) \to \bZ/3\bZ.
\]
The map $\ph$ is up to isomorphism the unique $\Gal_\bQ$-map $\rT_3 E \to \bZ/3\bZ$ according to Lemma \ref{lem:E3}.

\begin{lem} \label{lem:goodreduction}
The $\bZ/3\bZ$ cover $h:W \to S_k$ has good reduction outside primes dividing $2\cdot 3 \cdot 5$.
\end{lem}
\begin{pro}
As $K/k$ is unramified above the primes in question, we may replace $h$ by the base change $W\times_k K \to S \times_\bQ K$. Here, the cover is the normalisation in the extension of the function field given by $\sqrt[3]{F}$ which by (\ref{eq:defF}) amounts to the same as by $\sqrt[3]{ \bruch{\sigma^2(\gamma)}{\gamma} \cdot  \bruch{2Y + \ep X}{2Y + \zeta^2 \ep X}}$. We have to check for prime divisors with multiplicity in $\divisor(\bruch{\sigma^2(\gamma)}{\gamma} \cdot  \bruch{2Y + \ep X}{2Y + \zeta^2 \ep X})$ not divisible by $3$, hence only $\divisor(\sigma^2(\gamma)/\gamma)$, which lives only in the fibres above $2$ and $5$ by construction.
\end{pro}


\subsection{The global Brauer class}
The global Brauer class, that will occur in the  Brauer--Manin obstruction against rational points on the Selmer curve $S$ will be the image of the cup product $(a,\cores_{k/\bQ}(m_F)) \in \rH^2(\pi_1(S),\mu_3)$ for some suitable $a \in \bZ[1/30]^\ast$, more precisely the inflation of its associated Kummer character in $\rH^1(\bQ,\mu_3)$.
For an adelic section $(s_p) \in \cS_{\pi_1(S/\bQ)}(\bA_\bQ)$ we have the adjunction formula
\begin{equation} \label{eq:adj}
\langle (a,\cores_{k/\bQ}(m_F)) , (s_p) \rangle = \sum_p \inv_p \big(s_p^\ast  (a,\cores_{k/\bQ}(m_F)) \big)
= \sum_p \inv_p \big((a,s_p^\ast  \cores_{k/\bQ}(m_F)) \big)
\end{equation}
\[
 = \sum_p \inv_p \big((a, \cores_{k/\bQ}(m_F \circ s_p|_{\Gal_{k_v}})) \big)  = \sum_p \inv_p \big(\cores_{k/\bQ}(a, m_F \circ s_v) \big)
\]
\[
= \sum_v \inv_v \big((a,m_F \circ s_v) \big) = \langle (a,m_F) , (s_v) \rangle
\]
with $s_v = s_p|_{\Gal_{k_v}}$ for a place $v|p$ of $k$ with completion $k_v$. Consequently, instead of obstructing $\bQ$-rational points, or sections, via $(a,\cores_{k/\bQ}(m_F))$, we can and will discuss the obstruction against $k$-rational points or sections of $\pi_1(S_k/k)$ via the class $(a,m_F)$ for $a \in \bQ^\ast$. Because $S$ as a cubic has $\bQ$-rational divisors of degree $3$, the class of the extension $\pi_1(S/\bQ)$ has order dividing $3$ and splits if and only if $\pi_1(S_k/k)$ splits. 

Multiplying with our fixed $\zeta_3$ we can transform the local invariants $\inv_v\big((a,m_F \circ s_v) \big)$ into the local cubic Hilbert pairing $(a,F(s_v))$ with values in $\mu_3$. Here $F(s_v)$ is the abuse of notation for the class $s_v^\ast \chi_F = \chi_F \circ s_v \in \rH^1(k_v,\mu_3)$.
The local contribution for a section $s_{x_v}$ associated to a point $x_v \in S(k_v)$ where $F$ is invertible, computes then as $(a,s_{x_v}^\ast(F)) = (a,F(x_v))$.

The alert reader will already have noticed that the Brauer--Manin obstruction that we will exploit for the Selmer curve amounts to the descent obstruction imposed by $h: W \to S_k$. For more details on the descent obstruction and the section conjecture we refer to \cite{stix:evidence}.


\subsection{Variance of the evaluation at points versus at sections}

The computation of (\ref{eq:adj}) shows moreover that the value of the local components of the Brauer--Manin obstruction depend only on $m_F \circ s_v \in \rH^1(k_v,\bZ/3\bZ)$, which subsequently is paired with the class of $a$ via the local Tate-duality pairing
\begin{equation}
\rH^1(k_v,\bZ/3\bZ) \times \rH^1(k_v,\mu_3) \to \bQ/\bZ.
\end{equation}
Let us identify $E(k_v) \otimes \bZ_3$ via the Kummer sequence of $E$ with its image in $\rH^1(k_v,\rT_3 E)$. The range of $m_F \circ s_v$, where $s_v$ ranges over all sections associated to $k_v$-points of $S$ (resp.\ all sections) is given as an affine subspace
\begin{equation} \label{eq:vark}
m_F \circ s_0 + \ph\big(E(k_v) \otimes \bZ_3\big) \subseteq \left( \text{resp. } m_F \circ s_0 + \ph( \rH^1(k_v,\rT_3 E))\right) \subseteq \rH^1(k_v,\bZ/3\bZ),
\end{equation}
where $s_0$ is the section for some point in $S(\bQ_p)$ for $v|p$. As $\Gal(k/\bQ)$ is of order prime to $3$, the corresponding situation over $\bQ_p$ is just given by the  $^+$-eigenspaces for the Galois action:
\begin{equation} \label{eq:varQ}
m_F \circ s_0 + \ph\big(E(\bQ_p) \otimes \bZ_3\big) \subseteq \left( \text{resp. } m_F \circ s_0 + \ph( \rH^1(\bQ_p,\rT_3 E))\right) \subseteq \rH^1(\bQ_p,\bZ/3\bZ).
\end{equation}
For the latter it was essential to choose the section of reference $s_0$ associated to a $\bQ_p$-point of $S$.


\subsubsection{Explicit computation of the variance}
The isogeny $f: E' \to E$ yields a short exact sequence of $\Gal_\bQ$-modules
\begin{equation} \label{eq:tatemodules}
0 \to \rT_3 E' \xrightarrow{\pi_1(f)} \rT_3 E \xrightarrow{\ph} \bZ/3\bZ \to 0.
\end{equation}
Before we exploit this cohomologically we need several observations. After base change to $k$ the isogeny of degree $3$ over $E_k$ is still unique up to isomorphism by Lemma \ref{lem:E3}. But complex multiplication by $\bZ[\zeta_3]$ is now defined, so that $f$ becomes isomorphic to the multipliction map $[\zeta_3 -1] : E_k \to E_k$, which is of degree $\deg [\zeta_3-1] = N_{k/\bQ} (\zeta_3 -1) = 3$. A description over $\bQ$ is as follows.

\begin{lem} \label{lem:isog}
\begin{enumerate}
\item  The elliptic curve $E'$ is the  quadratic twist of $E$ by the cyclotomic character $\Gal_{\bQ} \to \Gal(k/\bQ)=\{\pm 1\}$.
\item In explicit equations we have $E' = \{v^2 = u^3+900\}$ and $E=\{b^2=a^3-24300\}$ with the isogeny $f:E' \to E$ being isomorphic to the map 
\begin{equation} \label{eq:isog}
(u,v) \mapsto (a,b) = (\bruch{u^3+3600}{u^2},v\cdot\bruch{u^3-7200}{u^3}).
\end{equation}
\item The map $f$ is surjective on $\bQ_3$-points.
\end{enumerate}
\end{lem}
\begin{pro}
(1) follows visibly from (2). The isomorphism of $E$ as a cubic $\{A^3+B^3+60C^3=0\}$ to the Weierstra\ss-form given in the proposition comes by the map
\[
[A:B:C] \mapsto (a,b) = (-180 \cdot \bruch{C}{A+B}, 270 \cdot \bruch{A-B}{A+B}).
\]
The map defined by (\ref{eq:isog}) is a degree $3$ isogeny, which we know is unique up to isomorphism.

In order to prove the surjectivity of $f$ on $\bQ_3$ points we may argue for all but finitely many points $(a,b) \in E(\bQ_3)$. We have to solve the equations 
\[
a=  \bruch{u^3+3600}{u^2} \quad \text{ and }  \quad 
b=  v\cdot\bruch{u^3-7200}{u^3}
\]
in $u,v \in \bQ_3$ which automatically (up to maybe the finitely many exceptions when $u^3=7200$) yields a point $(u,v) \in E'(\bQ_3)$. The equation for $v$ is linear, so we need to worry only  about the equation for $u$. Substituting $T = u/a$ (except for the case $a=0$) we need to guarantee a solution in $\bQ_3$ of the equation
\begin{equation} \label{eq:preimage}
T^3 - T^2 +3600/a^3 = 0,
\end{equation}
as long as $a$ is part of a point $(a,b) \in E(\bQ_3)$. As $v_3(b^2)$ is even and $v_3(24300) = 5$, the triangle equality holds with a minimum $v_3(a^3) = \min\{v_3(b^2),5\}$. But $v_3(a^3)$ is divisible by $3$, thus $v_3(a^3) = v(b^2) < 5$, and so the value is divisible by $6$ and negative.
Now equation (\ref{eq:preimage}) has coefficients in $\bZ_3$ and reduces to $T^3-T^2$ modulo $3$. The separable solution $T=1$ lifts to a solution of (\ref{eq:preimage}) by Hensel's Lemma, which proves the lemma.
\end{pro}

\begin{lem} \label{lem:cyclic}
Let $M \not= (0)$ be a $\bZ_3[\zeta_3]$-module of finite length, which is cyclic as an abelian group. Then $M$ is isomorphic to $\bZ/3\bZ$.
\end{lem}
\begin{pro}
We have $M \cong \bZ_3[\zeta_3]/(\zeta_3-1)^n$ for some $n$. If $n\geq 2$, then we have a quotient $M \surj \bZ_3[\zeta_3]/(\zeta_3-1)^2 =\bZ/3\bZ \oplus \bZ/3\bZ \cdot \zeta_3$ which fails to be cyclic as an abelian group. Hence $n=1$, which proves the lemma.
\end{pro}

\begin{prop} \label{prop:notacube}
Let $v$ be a place of $k$ such that $60$ is not a cube in $k_v$.
\begin{enumerate}
\item The map $\rH^2(k_v,\rT_3 E) \to \rH^2(k_v,\bZ/3\bZ)$ induced by $\ph$ is an isomorphism of $\Gal_\bQ$-modules isomorphic to $\mu_3$.
\item $\rH^2(k_v,\rT_3 E') \cong \bZ/3\bZ$ with trivial action of $\Gal(k/\bQ)$.
\item The sequence of $\Gal(k/\bQ)$-modules 
\[
\rH^1(k_v,\rT_3 E) \xrightarrow{\ph} \rH^1(k_v,\bZ/3\bZ) \xrightarrow{(-,60)} \bZ/3\bZ \to 0
\] 
is exact, where $(-,60)$ is the pairing with the Kummer character associated to $60$.
\item $\ph(\rH^1(k_v,\rT_3 E)) \subset \rH^1(k_v,\bZ/3\bZ)$ is equal to the annihilator of  $\langle 60 \rangle \subset \rH^1(k_v,\mu_3)$ under the local Tate-duality pairing.
\item If $v \nmid 3$, then $E(k_v) \otimes \bZ_3 = \rH^1(k_v,\rT_3 E)$.
\end{enumerate}
\end{prop}
\begin{pro}
(1) By local Tate-duality the dual map equals $\mu_3(k_v) \to E[3^\infty](k_v)$, which is injective. By Lemma \ref{lem:E3} and Lemma \ref{lem:cyclic} the $k_v$-rational $3$-primary torsion $E[3^\infty](k_v)$ is cyclic of order $3$. 

(2) This follows from (1) as $E'$ is the quadratic twist of $E$ by Lemma \ref{lem:isog}.

For (3) we use the cohomology sequence for (\ref{eq:tatemodules}) and (1) and (2), and it only remains to compute the boundary map $\rH^1(k_v,\bZ/3\bZ) \to \rH^2(k_v,\rT_3 E')$. The short exact sequence (\ref{eq:tatemodules}) maps canonically to (\ref{eq:E3}) inducing isomorphisms on the relevant terms of the cohomology sequence by (2), so that the boundary map equals
\[
- \cup 60 \ : \ \rH^1(k_v,\bZ/3\bZ) \to \rH^2(k_v,\mu_3),
\]
which is nothing but the pairing map with the Kummer class associated to $60$, the class of the extension (\ref{eq:E3}) following Lemma \ref{lem:E3}. (4) is merely a reformualtion of (3).

(5) The 'multiplication by $3^n$'-sequences and continuous cohomology \`a la Jannsen yield the  short exact sequence
\[
0 \to E(k_v) \otimes \bZ_3 \to \rH^1(k_v,\rT_3 E) \to \rT_3(\rH^1(k_v,E)) \to 0.
\]
By local Tate-duality $\rH^1(k_v,E)$  equals $\Hom(E(k_v),\bQ/\bZ)$, which has finite $3$-primary torsion part and thus $\rT_3(\rH^1(k_v,E)) = \Hom(\bQ_3/\bZ_3,\rH^1(k_v,E))$ vanishes.
\end{pro}

\begin{prop} \label{prop:above3}
Let $v = (\zeta_3-1)$ be the place of $k$ dividing $3$.
\begin{enumerate}
\item  $\ph(E(k_v) \otimes \bZ_3) \subset \rH^1(k_v,\bZ/3\bZ)$ is equal to the annihilator of  $\langle 2,3 \rangle \subset \rH^1(k_v,\mu_3)$ under the local Tate-duality pairing.
\item $\ph(E(\bQ_3) \otimes \bZ_3) = 0 $ in $\rH^1(\bQ_3,\bZ/3\bZ)$.
\end{enumerate}
\end{prop}
\begin{pro} The snake lemma applied to 
\begin{equation} \label{eq:Esnake}
\xymatrix@M+1ex@R-2ex{ 0 \ar[r] & E'(k_v) \otimes \bZ_3 \ar[d]^f \ar[r] & \rH^1(k_v, \rT_3 E') \ar[d]^f \ar[r] & \rT_3 (\rH^1(k_v,E')) \ar[d]^f \ar[r] & 0 \\
0 \ar[r] & E(k_v) \otimes \bZ_3 \ar[r] & \rH^1(k_v, \rT_3 E) \ar[r] & \rT_3 (\rH^1(k_v,E)) \ar[r] & 0}
\end{equation}
yields 
\[
 (E(k_v) \otimes \bZ_3)/f(E'(k_v) \otimes \bZ_3) = \ph(E(k_v) \otimes \bZ_3).
\]
As $\bZ_3[\zeta_3]$-module we have $E(k_v) \otimes \bZ_3 = \bZ_3[\zeta_3] \oplus \bF_3$ and $f$ corresponds to multiplication by $\zeta_3 -1$, at least after identifying $E'_k = E_k$ correctly, and hence  $\dim_{\bF_3} \ph(E(k_v) \otimes \bZ_3) = 2$.

The subgroup $\langle 2,3\rangle$ is the $^+$-part of $\rH^1(k_v,\mu_3)$ under the action of $\Gal(k/\bQ)$, so that its annihilator under the Tate-duality pairing equals the $^-$-part of $\rH^1(k_v,\bZ/3\bZ)$ and has dimension $2$ over $\bF_3$.  For (1) it therefore suffices to show that $\ph(E(k_v) \otimes \bZ_3)$ has trivial $^+$-eigenspace, which is exactly statement (2).

The analog of (\ref{eq:Esnake}) with $k_v$ replaced by $\bQ_3$ yields 
\[
 (E(\bQ_3) \otimes \bZ_3)/f(E'(\bQ_3) \otimes \bZ_3) = \ph(E(\bQ_3) \otimes \bZ_3),
\]
which vanishes by Proposition \ref{lem:isog} (3).
\end{pro}


\subsection{The locally constant Brauer class for the Selmer curve}

\begin{thm} \label{thm:russbselmer}
 The class $\alpha = (2,m_F) \in \rH^2(\pi_1 S_k , \mu_3)$ maps to a locally constant Brauer class $A$ in $\RussB(S_k)$, and the class $\cores_{k/\bQ}(\alpha) = (a,\cores_{k/\bQ}(m_F)) \in \rH^2(\pi_1 S , \mu_3)$ maps to a class in $\RussB(S)$ that generates $\RussB(S)/\Br(\bQ)$.
\end{thm}
\begin{pro}
Let $A_v$ be the pull back of $A$ to $\Br(S \times_\bQ k_v)$.
By Lichtenbaum duality and because $S$ is of genus $1$, the class $A_v$ is locally constant if and only if its evaluation at $k_v$ points is independent of the chosen point in $S(k_v)$.

The class $\alpha$ is unramified at places $v \nmid 2 \cdot 3 \cdot 5$ by Lemma \ref{lem:goodreduction}, so that $A_v(x_v)$ for $x_v \in S(k_v)$ vanishes. At $v \mid 5$ we find that $2$ is a cube and thus $\alpha_v = 0$. 

At $v \mid 2 $ we may apply Proposition \ref{prop:notacube}. Because $15$ is a cube in $\bQ_2$, the annihilator of $2$ and of $60$ under the local Tate-duality pairing agree, and so $A(x_2)$ is independent of $x_2 \in S(k_v)$. Let $w$ be a place above $v$ in $K$. In fact, the exact value is $0$, because $F(x_2)$ is a norm from $K_w/k_v$ and $10$ is a norm from $K/k$, which by local class field theory implies that 
\[
\langle 10 \rangle \cdot (k_v^\ast)^3 = N_{K_w/k_v} K_w^\ast \subset  k_v^\ast.
\]
But the cubic Hilbert symbol $(2,10)$  for $k_v$ vanishes, hence also $A(x_2)$. It remains to deal with the place $v$ above $3$. Here Proposition \ref{prop:above3}(1) tells us, that the evaluation of $A$ is independent of the point.

With the class $\alpha$ also $\cores_{k/\bQ}(\alpha)$ is locally constant. By  \cite{mazur:localtoglobal} Thm 1 and \S9, we know that $\Sha^1(\bQ,E) \cong \bZ/3\bZ \times \bZ/3\bZ$, hence $\RussB(S)/\Br(\bQ)$ is cyclic of order $3$, see Section \ref{sec:CTlco}. In order to prove the last statement of the theorem it thus suffices to show, that $\alpha$ does not come from $\Br(k)$. Let $\delta = \sqrt[3]{10}$ be a cube root of $10$ in $\bQ_3$. The value of $F$ at $[0:\delta:-2]$ can be effortlessly evaluated with SAGE as
\[
F([0:\delta:-2]) = N_{K/k}\left(\bruch{\delta^2 - \gamma \delta + \gamma \cdot \sigma(\gamma)}{\delta^2}\right) = -\bruch{9}{5} \zeta_3 \delta^2 + (\bruch{9}{5} \zeta_3 + \bruch{36}{5}) \delta -\bruch{ 81}{5} \zeta_3 + 9 
\]
which up to cubes in $\bQ_3(\zeta_3)$ equals $(\zeta_3 - 1)(1 + (\zeta_3 - 1)^2)$, and pairs nontrivially with $2$ under the local cubic Hilbert pairing of $\bQ_3(\zeta_3)$, which can be evaluated (using SAGE) via the explicit formula for the Hilbert pairing from \cite{fesenkovostokov} VII \S4.
\end{pro}

Our discussion of the arithmetic of the Selmer curve so far allows to imediately disprove any hope that sections of $\pi_1(S/\bQ)$ might be Brauer--Manin obstructed in the same way as rational points are. 

\begin{thm}
There is an adelic section $(s_p) \in \cS_{\pi_1(S/\bQ)} (\bA_\bQ)$ which survives the Brauer--Manin obstruction from classes of $\rH^2_{\RussB}(\pi_1(S),\mu_3)$, in particular from the class $\cores_{k/\bQ}(\alpha)$ of Theorem \ref{thm:russbselmer}.
\end{thm}
\begin{pro}
We have a short exact sequence
\[
0 \to \Pic(S) \otimes \bF_3 \to \rH^2_{\RussB}(\pi_1(S),\mu_3)/\rH^2(\bQ,\mu_3) \to \RussB(S)/\Br(\bQ) \to 0,
\]
in which $\Pic(S) \otimes \bF_3$ and $ \RussB(S)/\Br(\bQ)$ are both cyclic of order $3$. The Brauer--Manin obstruction is linear in the class from $\rH^2_{\RussB}(\pi_1(S),\mu_3)$ and vanishes on the image of $\rH^2(\bQ,\mu_3)$ by inflation. Hence it suffices to consider the Brauer--Manin obstruction imposed by the class $\cores_{k/\bQ}(\alpha)$ from Theorem \ref{thm:russbselmer}, and that imposed by a generator of $\Pic(S) \otimes \bF_3$.

By the proof of Theorem \ref{thm:russbselmer} and by Proposition \ref{prop:notacube} we find that the local components of $\langle \cores_{k/\bQ}(\alpha) , (s_p)\rangle$ vanish except for possibly the component above $3$. 

At $3$ Proposition \ref{prop:notacube}(4) shows, that the variation of $m_F \circ s_3$ with $s_3$ the $3$-adic component of an adelic section is the annihilator of $60$, which is a larger space than the variation of $m_F \circ s_{x_3}$ for $x_3 \in S(\bQ_3)$ after  Proposition \ref{prop:above3}. In particular, because 
\[
m_F \circ  s_{[0:\delta:-2]} \in \langle 60 \rangle^\perp = \ph(\rH^1(\bQ_3,\rT_3 E)) \oplus  \langle 2,3 \rangle^\perp,
\]
there is a $3$-adic section $s_3$, which necessary does not belong to a point in $S(\bQ_3)$, such that $m_F \circ s_3$ lies in the annihilator of $\langle 2,3 \rangle$. Hence, if we choose this section as the $3$-adic component, then the Brauer--Manin pairing with $\cores_{k/\bQ}(\alpha)$ vanishes.

Now we look at a generator of  $\Pic(S) \otimes \bF_3$, which we choose as the divisor of the homogeneous coordinate $X$, so for $\delta = \sqrt[3]{10}$ we look at 
\[
D = [0 : \delta : -2] + [0 : \zeta_3\delta : -2]  + [0 : \zeta_3^2 \delta : -2] .
\]
As $10$ becomes a cube in $\bQ_3$, we may assume $\delta \in \bQ_3$ and compute in $\Pic(S \otimes_\bQ \bQ_3) \otimes \bF_3$ 
\[
D \equiv D - 3 \cdot [0 : \delta : -2] = [0 : \zeta_3\delta : -2]  + [0 : \zeta_3^2 \delta : -2]  - 2\cdot  [0 : \delta : -2]  = \divisor(\bruch{X}{2Y + \delta Z}).
\]
Therefore the associated Brauer--Manin obstruction against adelic sections has no contribution at $3$. The summands at $p \not=3$ vanish as well, because by Proposition \ref{prop:notacube}(5) we may assume that the $p$-adic section belongs to a $p$-adic point, and for those vanishing is obvious by functoriality and $\Pic(\bQ) = 0$.
\end{pro}

The discussion of the Selmer curve in light of the section conjecture would not be complete without posing the following question, which the author so far is unable to decide.

\begin{ques}
Does $\pi_1(S/\bQ)$ split?
\end{ques}

Of course, the extension $\pi_1(S/\bQ)$ is the fibre product over $\Gal_\bQ$ of the geometrically pro-$3$ extension $\pi_1^3(S/\bQ)$ with the geometrically prime-to-$3$ extension $\pi_1^{3'}(S/\bQ)$. The latter splits due to a corestriction/restriction argument and by means of a point of $S$ in a field extension of $\bQ$ of degree $3$. Consequently, the splitting of $\pi_1(S/\bQ)$ is equivalent to the splitting of $\pi_1^3(S/\bQ)$.


\begin{appendix}


\section{Letter from Deligne to Thakur} \label{app:letterDeligne}

\noindent We reproduce the following letter by P.~Deligne to D.~Thakur with Deligne's gratefully acknowledged authorization.

\medskip

\hfill March 7, 2005

Dear Thakur,

\bigskip

I thought I had a proof, but it is wrong. It contained one observation which I still like, which follows.
\medskip

Let $X$ be smooth over a ring of $S$-integers $\dO$ of a number field $F$. Suppose $X_F/F$ absolutely irreducible. If we choose an algebraic closure $\bar{F}$ of $F$, and a geometric point $\bar{x}$ of $X_{\bar{F}}$, we get an exact sequence of profinite groups
\begin{equation} \label{eq:1}
1 \to \pi_1(\bar{X}) \to \pi_1(X) \to \pi_1(F) \to 1
\end{equation}
with $\pi_1(\bar{X}) := \pi_1(X_{\bar{F}},\bar{x})$, $\pi_1(X) := \pi_1(X,\bar{x})$ and $\pi_1(F) = \pi_1(\Spec(F),\Spec(\bar{F}))$. A point $s \in X(F)$ given with a path, on $X_{\bar{F}}$ from $s(\Spec(\bar{F}))$ to $\bar{x}$, defines a splitting of (\ref{eq:1}), and $s$ alone defines a $\pi_1(\bar{X})$-conjugacy class of splittings.
\medskip

\underline{Observation}: \textit{Let $S$ be the union of the conjugacy classes of splitting of (\ref{eq:1})
defined by points $s \in X(\dO)$. The set $S$ is equicontinuous.}
\medskip

In other words, for any open subgroup $K$ of $\pi_1(X)$, there exists an open subgroup $L$ of $\pi_1(F)$ such that all $s \in S$ define the same composite map $L \to \pi_1(F) \to \pi_1(X) \to \pi_1(X)/K$.

It follows from this observation that the closure $\bar{S}$ of $S$ (in the space of continuous maps from $\pi_1(F)$ to $\pi_1(X)$ is compact, and also that the quotient of $\bar{S}$ by $\pi_1(\bar{X})$-conjugation is compact. If $X/\dO$ is proper, we have $X(F) = X(\dO)$, and if all $\pi_1(\bar{X})$-conjugacy classes of sections correspond one to one to points, this will give $X(F)$ a compact (pro-finite) topology. I believe now that my conviction of long ago that this would give a contradiction if $X$ was a curve of genus $\geq 2$ and if $X(F)$ were (denumerably) infinite was wishful thinking.
\medskip

Here is the proof of the observation.

\bigskip

Shrink $\Spec(\dO)$ such that $X/\dO$ can be compactified with a relative divisor with normal crossing at infinity.
\medskip

The group $\pi_1(\bar{X})$ being finitely generated (as a profinite group), it has only a finite number of (open) subgroups of index $\leq n$. Let $K$ be their intersection. It is an open characteristic subgroup, and $\pi_1(\bar{X})^{(n)} := \pi_1(\bar{X})/K$ has an order all of whose prime divisors are $\leq n$. The group $\pi_1(\bar{X})$ is the projective limit of the $\pi_1(\bar{X})^{(n)}$. Dividing by $K$, we deduce from (\ref{eq:1}) an exact sequence
\begin{equation} \label{eq:2}
1 \to \pi_1(\bar{X})^{(n)} \to \pi_1(X)^{(n)} \to \pi_1(F) \to 1
\end{equation}
A splitting of (\ref{eq:1}) induces a splitting of (\ref{eq:2}), and we will show that the splittings in $S$ induce morphisms $\pi_1(F) \to \pi_1(X)^{(n)}$ which are all equal on a suitable subgroup of finite index of $\pi_1(F)$.
\medskip

Define $\dO^{(n)} = \dO[1/n!]$ and $\pi_1(\dO^{(n)}) = \pi_1(\Spec(\dO^{(n)}), \Spec(\bar{F}))$. The sequence (\ref{eq:2}) is the pull back by $\pi_1(F) \to \pi_1(\dO^{(n)})$ of a sequence
\begin{equation} \label{eq:3}
1 \to \pi_1(\bar{X})^{(n)} \to \pi_1(X_{\dO^{(n)}})^{(n)} \to \pi_1(\dO^{(n)}) \to 1
\end{equation}
where $\pi_1(\dO^{(n)})$ is the quotient of $\pi_1(X_{\dO^{(n)}}, \bar{x})$ by the image of $K$. A point $s \in X(\dO)$ (or just $X(\dO^{(n)})$ defines a $\pi_1(\bar{X})^{(n)}$-conjugacy class of splittings of (\ref{eq:3}), and the splitting of (\ref{eq:2}) deduced from the splitting of (\ref{eq:1}) in $S$ are pull back by $\pi_1(F) \to \pi_1(\dO^{(n)})$. It remains to prove the

\begin{lem*}
The set of splittings of the short exact sequence (\ref{eq:3}) is finite.
\end{lem*}
\begin{pro*}
By Minkowski, $\pi_1(\dO^{(n)})$ has only a finite number of (open) subgroup of a given index. Let $a$ be the order of $\Aut(\pi_1(\bar{X})^{(n)})$ and $z$ the order of the center of $\pi_1(\bar{X})^{(n)}$. Let $H$ be the intersection of the subgroups of $\pi_1(\dO^{(n)})$ of index dividing $a^2z$. It is a subgroup of finite index, and I claim that any two sections agree on $H$, a claim from which the lemma follows. Indeed, a section $s$ defines $a_s : \pi_1(\dO^{(n)}) \to \Aut(\pi_1(\bar{X})^{(n)})$. If $s$ and $t$ are sections, on $\ker(a_s) \cap \ker(a_t)$, $s$ and $t$ are with value in the centralizer of $\pi_1(\bar{X}^{(n)})$, and differ by an homomorphism to the center of $\pi_1(\bar{X})^{(n)}$. The section $s$ and $t$ hence agree on a subgroup of index dividing $a^2z$, and a fortiori on $H$.
\end{pro*}
\medskip

\hspace{0.5cm} All the best,
\medskip

\hspace{1cm} P.~Deligne

\end{appendix}







\end{document}